\documentclass[12pt,a4paper]{amsart}
\usepackage[top=1.15in, bottom=1.15in, left=1.17in, right=1.17in]{geometry}
\usepackage{amssymb}
\usepackage{pdflscape}
\usepackage{amscd}
\usepackage{fancyhdr}
\input xy 
\xyoption{all}

\usepackage{amsmath,amsfonts,stmaryrd}
\usepackage[usenames,dvipsnames,svgnames]{xcolor}
\usepackage[colorlinks, urlcolor=Navy, linkcolor=Navy, citecolor=Navy]{hyperref}
\usepackage[abbrev,msc-links]{amsrefs}
\usepackage[latin1]{inputenc}
\usepackage{tikz}
\usetikzlibrary{shapes,arrows}

\theoremstyle{plain}

\newtheorem{thm}{Theorem}[section]
\newtheorem{lem}[thm]{Lemma}

\newtheorem{pro}[thm]{Proposition}
\newtheorem{cor}[thm]{Corollary}
\newtheorem{conj}[thm]{Conjecture}
\theoremstyle{definition}
\newtheorem{dfn}[thm]{Definition}
\newtheorem{exa}[thm]{Example}
\newtheorem{rem}[thm]{Remark}
\newtheorem{que}[thm]{Open question}

\DeclareMathOperator{\Hom}{Hom}

\DeclareMathOperator{\End}{End}
\DeclareMathOperator{\Ext}{Ext}
\DeclareMathOperator{\Tor}{Tor}

\DeclareMathOperator{\modu}{mod}

\DeclareMathOperator{\Lamod}{\Lambda-\rm{mod}}

\DeclareMathOperator{\Modu}{Mod}

\DeclareMathOperator{\ModmcA}{\rm{Mod}-\mathcal{A}}

\newcommand{\op}{\mathrm{op}}

\DeclareMathOperator{\idim}{id}
\DeclareMathOperator{\pdim}{pd}
\DeclareMathOperator{\gldim}{gldim}

\newcommand{\oTo}{\xymatrix{ \ar@{^{(}->}[r]|{\mathbf{O}}& }} 
\newcommand{\cTo}{\xymatrix{ \ar@{^{(}->}[r]|{\mathbf{|}}& }} 
\newcommand{\coTo}{\xymatrix{ \ar@{^{(}->}[r]|{\mathbf{O}}|{\mathbf{|}}& }} 

\DeclareMathOperator{\Ker}{ker}
\DeclareMathOperator{\coKer}{coker}
\DeclareMathOperator{\Bild}{Im}

\DeclareMathOperator{\kdual}{D}

\DeclareMathOperator{\cogen}{cogen}
\DeclareMathOperator{\gen}{gen}

\newcommand{\La}{\Lambda}

\newcommand{\mcA}{\mathcal{A}}

\newcommand{\mcC}{\mathcal{C}}

\newcommand{\mcE}{\mathcal{E}}

\newcommand{\mcI}{\mathcal{I}}

\newcommand{\mcL}{\mathcal{L}}
\newcommand{\mcM}{\mathcal{M}}

\newcommand{\mcP}{\mathcal{P}}
\newcommand{\mcQ}{\mathcal{Q}}

\newcommand{\mcS}{\mathcal{S}}
\newcommand{\mcT}{\mathcal{T}}

\newcommand{\mcX}{\mathcal{X}}
\newcommand{\mcY}{\mathcal{Y}}

\DeclareMathOperator{\copres}{copres}
\DeclareMathOperator{\pres}{pres}
\DeclareMathOperator{\Cores}{Cores}
\DeclareMathOperator{\Reso}{Res}

\DeclareMathOperator{\Thick}{Thick}

\DeclareMathOperator{\Add}{Add}
\DeclareMathOperator{\add}{add}

\DeclareMathOperator{\tilt}{tilt}
\DeclareMathOperator{\Tilt}{Tilt}

\usepackage{todonotes}

\begin{document}

\title{Tilting theory in exact categories}
\author{Julia Sauter}
\address{Julia Sauter\\ Faculty of Mathematics \\
Bielefeld University \\
PO Box 100 131\\
D-33501 Bielefeld }
\email{jsauter@math.uni-bielefeld.de}

\begin{abstract}
 We define tilting subcategories in arbitrary exact categories to archieve the following. Firstly: Unify existing definitions of tilting subcategories to arbitrary exact categories. Discuss standard results for tilting subcategories: Auslander correspondence, Bazzoni description of the perpendicular category.\\
 Secondly: We treat the question of induced derived equivalences separately - given a tilting subcategory $\mcT$, we ask if a functor on the perpendicular category induces a derived equivalence to a (certain) functor category $\modu_{\infty}-\mcT$ over $\mcT$. If this is the case, we call the tilting subcategory \emph{ideq} tilting. 
    We prove a generalization of Miyashita's theorem (which is itself a generalization of a well-known theorem of Brenner-Butler) and characterize exact categories with enough projectives allowing ideq tilting subcategories.\\
    In particular, this is always fulfilled if the exact category is abelian with enough projectives. 
\end{abstract}


\maketitle

\section{Introduction}
Tilting theory ($=$ categories with tilting objects) is originally defined for categories of finitely generated modules over artin algebras (Brenner-Butler), here the tilting modules were still assumed of $\pdim \leq 1$. Then this was generalized to arbitrary projective dimension by \cite{Hap}, chapter 3 - here you also find a detailed account of the beginning of tilting theory (starting with BGP reflection functors). 
Afterwards this was generalized in several directions (  
\cite{W}, 
\cite{Mi}, \cite{CF}, 
\cite{AC}, \cite{CT}, 
\cite{HRS},
\cite{AHK} chapter 5, \cite{Ric} - referred to (by us) as: \emph{infinite} or \emph{big} tilting or tilting in triangulated categories). 
The different developments of the time (2007) were captured in the handbook of tilting theory \cite{AHK}. Since then many more generalizations were found, e.g. the later discussed recent works  \cite{ZZ}, \cite{AM}, \cite{Rump}.\\
For exact categories: The first occurrence in \cite{ASoII} is tilting objects in exact substructures of $\modu-\La$ for an artin algebra $\La$. 
In \cite{Kr-book}, chapter 7, a tilting object in an exact category is defined as a self-orthogonal object which generates the bounded derived category of the exact category as a thick subcategory. 
 In \cite{ZZ}, the authors define tilting subcategories for extriangulated categories with enough projectives and enough injectives. 
 The literature on \emph{infinite} or \emph{big} tilting in exact category will not be considered here (this includes \cite{Rump}). An alternative definition of tilting in exact categories can be found in \cite{AM} which covers big and small tilting in exact categories (as far as we can see: There is an additional axiom (T3) required in loc. cit. which we are not using). 

Our motivation is to generalize and unify the following classical and recent definitions and results (which we only sketch as follows): \\
\begin{itemize}
    \item[(A)] Tilting modules over artin algebras induce derived equivalences (\cite{Hap} chapter 3, \cite{Ric})
    \item[(B)] Relative tilting modules over artin algebras induce derived equivalences (\cite{ASoII}, \cite{Bu-Closed}) 
    \item[(C)] Tilting bundles over projective varieties over a field induce derived equivalences (\cite{Bae}, \cite{Bei})
    \item[(D)] Zhang-Zhu \cite{ZZ} introduce tilting subcategories for extriangulated categories with enough projectives and enough injectives. They prove a generalization of a result called Auslander-correspondence.  
\end{itemize}
To include (C),(D): We must drop the assumption that the exact category has enough projectives and we have to generalize from tilting modules to tilting subcategories. 
Our definition can be found as Def. \ref{tiltingDef}. Furthermore, we shortly discuss standard results from tilting theory: \\
1. \emph{Bazzoni's description of the perpendicular category} cf. Corollary \ref{Baz}. This gives a description of the perpendicular category of a tilting subcategory which is practical to find examples (see e.g. special tilting). \\
2. \emph{Auslander correspondence} cf. subsection \ref{Auslander-corresp}.  This means a characterization of the subcategories which arise as perpendicular categories of tilting subcategories. \\
As this topic has so many precursors, we discuss compatibility with other definitions of tilting (this does not claim completeness due to the amount of literature on the subject) in section \ref{DefOfTilting}.\\
To generalize (A),(B): We need to find a \emph{tilting functor} which induces in (A) and (B) a triangle equivalence as claimed.  \\
We introduce in section \ref{Ideq} the notion ideq tilting, which means that the tilting functor induces a triangle equivalence on the bounded derived categories. We prove, if for a an $n$-tilting subcategory $\mcT$, the category $\modu_{\infty}-\mcT$ has finite global dimension, then $\mcT$ is ideq tilting (cf. Prop. \ref{FinDerEq})- this generalizes (C). We prove the following results as generalizations of (A) and (B) respectively: 

\begin{thm} (cf. Thm \ref{ideqWep})
Let $\mcE$ be an exact category with enough projectives $\mcP$. Then the following are equivalent: 
\begin{itemize}
\item[(1)] $\mcE$ is equivalent as an exact category to a finitely resolving subcategory $\overline{\mcE}$ of $\modu_{\infty}-\mcP$, i.e. $\overline{\mcE}$ is resolving and for every $F\in \modu_{\infty}-\mcP$ there is a finite exact sequence $0 \to E_n \to \cdots \to E_0 \to F \to 0$ with some $E_i \in \overline{\mcE}$ for some $n \geq 0$.
\item[(2)] There is an $n\in \mathbb{N}_0$ and an $n$-tilting subcategory of $\mcE$ which is ideq $n$-tilting.
\item[(3)] For every $n\geq 0$, every $n$-tilting subcategory of $\mcE$ is ideq $n$-tilting. 
\end{itemize}
\end{thm}

\begin{cor} (cf. Cor. \ref{relIdeqWep})
  Let $\mcP$ be an idempotent complete, 
  additive category. 
 Let $\mcE$ be an exact substructure of $\modu_{\infty}-\mcP$, with enough projectives $\mcQ$. Then for every $n\geq 0$, every $n$-tilting subcategory of $\mcE$ is ideq $n$-tilting. 
 \end{cor}
 
To prove the first theorem, we prove a \emph{Miyashita Theorem} (generalization of Brenner Butler's theorem) cf. Theorem \ref{MiThm}. This describes the image of the perpendicular category of an ideq tilting subcategory $\mcT$ under the tilting functor ($X \to \Hom_{\mcE}(-,X)\lvert_{\mcT}$).

As standing assumption: We will always assume that the exact category is idempotent complete. 
\\
To introduce definitions of subcategories and recall results from the literature, we start in section 2 with preliminaries on subcategories of exact categories and in section 3 we give a quick introduction to the bounded derived category. \\
The author is supported by the Humboldt Professorship of William Crawley-Boevey and would like to thank him for helpful discussions.

\section{Some definitions of subcategories}

Let $\mcA$ be an idempotent complete category. 
Let $\mcE=(\mcA, \mcS)$ be throughout this section be an exact category in the sense of Quillen (it consists of an additive category $\mcA$ together with a class of kernel-cokernel pairs $\mcS$, referred to as \textbf{short exact sequences} which satisfy the axioms of \cite{Bue}, Def. 2.1).\\
For a kernel-cokernel pair $(i,d)\in \mcS$ we call $i$ an \textbf{inflation} and $d$ a \textbf{deflation}.
We will denote by $\mcP(\mcE)$ the projective objects in $\mcE$ and by $\mcI(\mcE)$ the injectives. Since this is common practice, we will also often denote the underlying additive category $\mcA$ again by $\mcE$ - we think the reader can handle this level of ambiguity. 

\begin{dfn}
If $\mcE=(\mcA, \mcS)$ is an exact category and $\mcX$ a full subcategory which is closed under extensions. Then we call $\mcX$ \textbf{fully exact subcategory} if we consider it together with the exact structure $\mcS\lvert_{\mcX}$ (i.e. the short exact sequences in $\mcS$ where all three terms lie in $\mcX$). We will write $\mcP(\mcX)$ for the Ext-projectives in $\mcX. $ 
\end{dfn}

\subsection{Subcategories generated by or orthogonal to a subcategory}

We call a morphism $f\colon X \to Y$ in $\mcE$ \textbf{admissible} if it factors as $f= d\circ i$ for an inflation $i$ and a deflation $d$. We say $X\xrightarrow{f} Y \xrightarrow{g}Z$ with $f,g$ admissible is exact at $Y$ if $\Bild f = \Ker g$ and $\Ker g \to Y \to \Bild g$ is an exact sequence in $\mcE$. A sequence of composable morphisms is \textbf{exact} if every morphism is admissible and the sequence is exact at every intermediate object (for short exact sequences we sometimes leave out the zeros in the beginning and end). We call a sequence 
\[ X_n \xrightarrow{f_n} X_{n-1} \xrightarrow{f_{n-1}} \cdots \to X_0 \xrightarrow{f_0} X_{-1}\xrightarrow{f_{-1}} 0\]
\textbf{right exact} if there is an exact sequence 
\[ 0 \to Z \to X_n \xrightarrow{f_n} X_{n-1} \xrightarrow{f_{n-1}} \cdots \to X_0 \xrightarrow{f_0} X_{-1}\xrightarrow{f_{-1}} 0\]
We say a (co- or contravariant) functor $F\colon \mcE \to (Ab)$ into abelian groups \textbf{is exact on} the right exact sequence if $F$ maps all short exact sequences 
 \[
 \Ker f_i \to X_i \to \Bild f_i, \quad -1 \leq i \leq n
 \]
 to short exact sequences in abelian groups. \\
Let $\mcX$ be a full additive subcategory of $\mcE$ and an integer $n\geq 0$ and a subset $I\subset \mathbb{N}_0$, we define the following full subcategories of $\mcE$ 
\[ 
\begin{aligned}
\gen_n(\mcX) &= \{ M \in \mcE \mid \exists \text{ right exact } X_n \to \cdots \to X_0 \to M \to 0,  \; X_i \in \mcX \\
 &{}\quad {} \quad {} \quad {}\quad {}\quad {} \quad \Hom_{\mcE}(X, -) \text{ exact on it for every }X \in \mcX \}\\
\pres_n (\mcX) &= \{ M \in \mcE \mid \exists \text{ right exact } X_n \to \cdots \to X_0 \to M \to 0,  \; X_i \in \mcX \} \\
\Reso_n (\mcX ) &= \{ M \in \mcE \mid \exists \text{ exact }0\to  X_{n} \to \cdots \to X_0 \to M \to 0,  \; X_i \in \mcX \}  \\
\mcX^{\perp_I} &=\{ M \in \mcE \mid \Ext_{\mcE}^i (X, M)=0 \text{ for all }i \in I, X \in \mcX \}
\end{aligned}
\]
We write $\mcX^{\perp_{\geq n}}:=\mcX^{\perp_{\left[ n,\infty\right)}}, \mcX^{\perp}:= \mcX^{\perp_{\geq 1}}$. We define $\gen_{\infty}(\mcX)$ (resp. $\pres_{\infty}(\mcX)$) analogously to $\gen_n(\mcX)$ (resp. $\pres_n(\mcX)$) with infinite resolutions. \\
We define $\gen (\mcX):= \gen_0 (\mcX), \pres (\mcX):= \pres_0 (\mcX)$ and observe $\gen (\mcX)= \pres (\mcX)$ is equivalent to $\mcX$ is contravariantly finite in $\pres (\mcX)$. Also, we set 
\[  
\Reso (\mcX) := \bigcup_{n\geq 1} \Reso_n(\mcX) 
\]
We have obvious inclusions 
\[ 
\Reso_n (\mcX ) \subset \pres_n (\mcX) \supset \gen_n (\mcX),\quad  \Reso (\mcX ) \subset \pres_{\infty} (\mcX) \supset \gen_{\infty} (\mcX)
\]
We denote the dual notions with $\cogen^n (\mcX), \copres^n(\mcX), n\leq \infty$ and $\Cores_n (\mcX), \Cores(\mcX), {}^{{}_I\perp}\mcX$ respectively.\\
Observe that these categories depend on the exact structure $\mcS$, so if there is the possibility of confusion, we will endow them with an index $\mcE$ (or $\mcS$).

\begin{dfn}
Let $X$ be an object in $\mcE$. We say $\pdim_{\mcE} X \leq n$ (resp. $\idim_{\mcE} X\leq n$) if $\Ext_{\mcE}^{n+1}(X,-)=0$ (resp. $\Ext_{\mcE}^{n+1} (-,X)=0$). \\
For a subcategory $\mcX$ we define 
\[ 
\pdim_{\mcE} \mcX:= \sup  \{\pdim_{\mcE} X  \mid X\in \mcX\} \quad \in \mathbb{N}_0 \cup \{ \infty\}
\]
and analogously $\idim_{\mcE} \mcX$. We call 
$\mcP^{\leq n}$ the full subcategory of all objects $X$ with $\pdim_{\mcE}X\leq n$ and $\mcP^{<\infty}$ the subcategory of all objects $X$ with $\pdim_{\mcE} X<\infty$. The subcategories $\mcI^{\leq n}$, $\mcI^{<\infty}$ are defined dually. 
\end{dfn}

\begin{rem}
We would like to know the common definition of an exact category $\mcE$ with "$\mcE = \mcP^{<\infty}$" since we use it so frequently. \\
If the exact category has the Jordan-H\"older property and only finitely many $\mcE$-simples, then $\gldim \mcE = \max_{S \text{ simple}} \pdim S$. Therefore, in that case $\mcE=\mcP^{<\infty}$ is equivalent to $\gldim \mcE < \infty$.\\
The \emph{finitistic dimension conjecture} holds for $\mcE$ if \[ \gldim (\mcP^{<\infty}) <\infty.\] So, also in this case, $\mcE=\mcP^{<\infty}$ is equivalent to $\gldim \mcE <\infty$.\\
(As far, as we know this conjecture is open for $\mcE=\modu-\La$ where $\La$ is an artin algebra. But is known to fail in some other abelian categories.) 
\end{rem}


We observe the following lemma.
\begin{lem} \label{basicProp}
Let $\mcE$ be an exact category and $\mcX$ a full, self-orthogonal subcategory. 
\begin{itemize}
    \item[(a)] Then we have $\Reso (\mcX) \subset \mcX^{\perp}$ and $\Reso (\mcX )$ is an extension-closed subcategory. 
    \item[(b)] Let $n\geq 1$. If $\pdim_{\mcE} \mcX \leq n$, then we have $\pres_{n-1}(\mcX) \subset \mcX^{\perp}$ and $\pres_m(\mcX)$ is extension-closed for all $n-1 \leq m \leq \infty $. 
\end{itemize}
\end{lem}

\begin{proof}
ad (a) Since $\mcX^{\perp}$ is inflation-closed it follows $\Reso (\mcX) \subset \mcX^{\perp}$. To see that $\Reso (\mcX) $ is extension closed one can use literally proof of the horseshoe lemma (replacing the projectives with $\mcX$), cf. \cite{Bue}, Thm. 12.8. \\
ad (b) Given an exact sequence $0\to Z \to X_{n-1} \to \cdots \to X_0 \to S \to 0$ with $X_i \in \mcX$, we see by dimension shift that $S \in \mcX^{\perp}$ using $\pdim_{\mcE} \mcX \leq n$. This proves $\pres_{n-1} (\mcX) \subset \mcX^{\perp}$ and therefore $\pres_m (\mcX) \subset \mcX^{\perp}$, $n-1\leq m \leq \infty$. Again by the horseshoe lemma as in loc. cit it follows that $\pres_m(\mcX)$ is extension closed. \\
\end{proof}

\subsection{Thick subcategories}
Recall that a full subcategory $\mcL$ of $\mcE$ is called \textbf{thick}, if it is closed under summands and for every short exact sequence $A \to B \to C $ in $\mcS$ the following holds: If any two objects of $A,B,C$ are in $\mcL$, then the third is also in $\mcL$.\\
We denote by $\Thick (\mcX)$ the smallest thick subcategory of $\mcE$ that contains $\mcX$. By definition, we have $\Reso (\mcX)\subset \Thick (\mcX)$ and $\Cores (\mcX) \subset \Thick (\mcX)$.  

\begin{exa}
The categories $\mcP^{<\infty}$ and $\mcI^{<\infty}$ are thick subcategories of $\mcE$. 
\end{exa}

\begin{dfn}
Let $\mcE$ be an exact category and $\mcX$ a full subcategory. We call $\mcX$ \textbf{inflation-closed} if for every inflation $i\colon X\to Y$ with $X,Y$ in $\mcX$ one also has $\coKer i $ in $\mcX$. Dually, one defines \textbf{deflation-closed}. 
\end{dfn}

\begin{pro} (\cite{AB}, Prop. 3.5 and Prop. 3.6, Thm 1.1) \label{AB-thick}
Let $\mcE$ be an exact category and $\mcX$ a fully exact subcategory which has enough projectives $\mcP=\mcP(\mcX)$. \\
Assume $\mcX$ is an inflation-closed subcategory. 
Then we have 
\begin{itemize}
\item[(1)] $\Thick (\mcX) =\Cores (\mcX)$, 
\item[(2)] $\Cores_n(\mcP)= \Cores_n(\mcX)\cap {}^{\perp} \mcX$ and therefore also $\Cores (\mcP)=\Cores (\mcX)\cap {}^{\perp}\mcX$, 
\item[(3)] $\mcX$ is covariantly finite in $\Cores (\mcX)$
\end{itemize}
\end{pro}

\begin{exa}
Let $\mcE$ be an exact category. We denote by $\mcP=\mcP(\mcE )$ the projectives. Then, $\mcP$ fulfills the dual conditions of Prop. \ref{AB-thick} and therefore 
\[ \Thick(\mcP)=\Reso(\mcP) \]
are all objects which admit a finite projective resolution. If $\mcE$ has enough projectives, then $\mcP^{<\infty}=\Reso (\mcP)$. 
\end{exa}


\subsubsection{Thick subcategories in triangulated categories}
Let $\mcC$ be a triangulated category. We call a full subcategory \textbf{thick} if it is closed under summands and a triangulated subcategory. \\
Let $\mcX$ be an additive subcategory of $\mcC$, we write $\Thick_{\Delta} (\mcX)$ for the smallest thick subcategory of $\mcC$ containing $\mcX$. 

\subsection{Resolving subcategories}

\begin{dfn} Let $\mcE$ be an exact category. 
Let $\mcX$ be a full subcategory. \\
We say $\mcX$ is \textbf{resolving} if it is extension closed, deflation-closed and $\pres (\mcX) = \mcE$.\\
We say $\mcX$ is \textbf{coresolving} if it is extension closed, inflation-closed and $\copres (\mcX) = \mcE$.\\
We say $\mcX$ is 
\begin{itemize}
    \item[(*)] \textbf{finitely resolving} if it is resolving and $\Reso (\mcX) = \mcE$,
    \item[(*)] $n-$\textbf{resolving} if it is resolving and $\Reso_n (\mcX) = \mcE$,
    \item[(*)] \textbf{uniformly finitely resolving} if it is $n$-resolving for some $n\geq 0$.
    \end{itemize}
Dually one defines finitely coresolving, $n$-coresolving and uniformly finitely coresolving. 
\end{dfn}

\begin{exa}
Let $\mcE$ be an exact category with enough projectives and $\mcX$ be any additive subcategory. Then ${}^{\perp}\mcX$ is a resolving subcategory. 
Dually, if $\mcE$ has enough injectives and $\mcX$ is an additive subcategory, then $\mcX^{\perp}$ is a coresolving subcategory. 
\end{exa}

\subsection{Functor categories}
Let $\mcX$ be an additive category. 
If $\mcX$ is essentially small we can define the category $\Modu-\mcX$ of right $\mcX$-modules as the category of contravariant additive functors on $\mcX$ to abelian groups (if $\mcX$ is not essentially small the class of all natural transformations between two $\mcX$-modules is not necessarily a set). This is an abelian category with enough projectives and injectives. The projectives are summands of arbitrary direct sums of representable functors $\Hom_{\mcX}(-,X)$ for some $X \in \mcX$ which we call the \emph{finitely generated projectives}. 
Similarly, one can define the category $\mcX-\Modu:= \Modu-\mcX^{op}$ of left $\mcX$-modules.  We define full  subcategories 
\[ 
\modu_{\infty}-\mcX \subset \modu_n-\mcX \subset \Modu-\mcX, \quad n \in \mathbb{N}_0
\]
via $F \in \modu_n-\mcX$ if there is a projective $n$-presentation (indexed by $0$ up to $n$) by finitely generated projectives, i.e. 
there is an exact sequence 
\[ 
\Hom_{\mcX}(-,X_n) \to \cdots \to \Hom_{\mcX} (-, X_0) \to F\to 0,
\]
we define $\modu_{\infty} -\mcX$ analogue with an infinite sequence as above and we call $\modu-\mcX:=\modu_0-\mcX$ the \emph{finitely generated} $\mcX$-modules, $\modu_1-\mcX$ as the \emph{finitely presented} $\mcX$-modules. All these are fully exact subcategories of $\Modu-\mcX$. \\

\begin{dfn}
Let $f\colon X_1 \to X_0$ be a morphism in $\mcX$, we say that $f$ \textbf{admits a weak kernel} if there exists a morphism $g \colon X_2 \to X_1$ such that 
\[ 
\Hom_{\mcX} (-,X_2 ) \xrightarrow{g\circ }\Hom_{\mcX} (-,X_1) \xrightarrow{f\circ } \Hom_{\mcX} (-,X_0)
\]
is exact in $\Modu-\mcX$. We say that $\mcX$ \textbf{admits weak kernels} if every morphism in $\mcX$ admits one. 
\end{dfn}

\begin{thm} (\cite{Au}, Prop. 2.1) The category $\modu_1-\mcX$ equals $\modu_{\infty}-\mcX$ if and only if $\mcX$ admits weak kernels. In this case $\modu_{\infty}-\mcX$ is abelian. 
\end{thm}

\begin{rem}
It is common practice to ignore the assumption $\mcX$ \emph{essentially small} because as soon as one looks at $\modu_1-\mcX$ the class of all natural transformations between two finitely presented $\mcX$-modules is a set. 
\end{rem}

\subsection{Embedding exact categories with enough projectives into a functor category}

If $\mcE$ is an exact category with enough projectives $\mcP$, then we the Yoneda embedding induces by composition with the restriction-on-$\mcP$ functor a functor 
\[ 
\mathbb{P} \colon \mcE \to \modu_{\infty}-\mcP, \quad E \mapsto \Hom_{\mcE}(-,E)\lvert_{\mcP}
\]
We write $\Bild \mathbb{P}$ for the essential image. For us, the main observation is the following:
\begin{pro} (\cite{E-master} Prop. 2.2.1, Prop. 2.2.8) 
The functor $\mathbb{P}$ is fully faithful and induces isomorphisms on all higher extension groups, $\Bild \mathbb{P}$ is extension closed. \\
If $\mcE$ is idempotent complete, then $\Bild \mathbb{P}$ is a resolving subcategory of $\modu_{\infty}-\mcP$.
\end{pro}

We observe also the following obvious lemma. 
\begin{lem} 
Given two exact categories $\mcE$ and $\mcE'$ with enough projectives $\mcP$ and $\mcP'$ respectively. Then the following are equivalent: 
\begin{itemize}
    \item[(1)] $\mcE$ and $\mcE'$ are equivalent as exact categories. 
    \item[(2)] There is an equivalence $\Psi\colon \mcP \to \mcP'$ of additive categories with the property: A morphism $f$ in $\mcP$ is admissible in $\mcE$ if and only if $\Psi (f)$ in $\mcP'$ is admissible in $\mcE'$. 
    \item[(3)] There is an additive equivalence $\Phi \colon \mcE \to \mcE'$ with $\Phi (\mcP) = \mcP'$.
\end{itemize}
\end{lem}

\begin{proof}
(2) implies (1): 
Given $\Psi$ as in (2). We extend $\Psi$ to an equivalence $\Psi \colon \mcE \to \mcE'$. For an (in $\mcE$) admissible $f\colon P_1\to P_0$ we map $\coKer (f)$ to $\coKer \Psi (f)$ (and then check that this is well-defined and gives an equivalence of categories). \\
Conversely, if $\Psi\colon \mcE\to \mcE'$ is an equivalence of exact categories it maps admissible morphsims to admissible morphisms, projectives to projectives. Then it restricts to an equivalence as in (2).\\
(3) implies (1): An additive equivalence preserves kernel-cokernel pairs. Kernel-cokernel pairs are short exact sequences in $\mcE$ (resp. $\mcE'$) if and only if $\Hom (P,-)$ are exact on them for all $P \in \mcP$ (resp. $\mcP'$). Therefore $\Phi$ is an equivalence of exact categories. 
\end{proof}

\begin{cor}
Let $\mcE$ be an exact category with enough projectives $\mcP$. 
The following are equivalent: 
\begin{itemize}
    \item[(a)] $\mathbb{P}\colon \mcE \to \modu_{\infty}-\mcP$ is an equivalence. 
    \item[(b)] The $\mcE$-admissible morphisms in $\mcP$ are precisely those such that $\Ker \Hom_{\mcP}(-,h) (\in \Modu -\mcP)$ lies in $\modu_{\infty}-\mcP$. 
\end{itemize}
\end{cor}

\begin{proof}
Clearly, we have (b) implies (a) and (a) implies (b) by the previous Lemma. 
\end{proof}

\begin{lem} \label{finResolving}
Let $\mcE$ be an exact category with enough projectives $\mcP$. Then the following are equivalent:  
\begin{itemize}
    \item[(1)] $\Bild \mathbb{P}$ is finitely resolving in $\modu_{\infty}-\mcP$.
    \item[(2)] For every morphism $f\colon P_1\to P_0$ in $\mcP$ admitting an infinite sequence of successive weak kernels in $\mcP$ (i.e. $\Hom_{\mcP}(-,f)$ is part of a projective resolution by finitely generated projectives of an object in $\modu_{\infty}-\mcP$) there exists a complex in $\mcE$
    \[ 
    0\to M_n \to M_{n-1} \to \cdots \to M_2 \to P_1 \xrightarrow{f} P_0 
    \]
    depending on $f$ 
    such that 
    \[ 
    0\to \mathbb{P} (M_n) \to \mathbb{P} (M_{n-1}) \to \cdots \to \mathbb{P}(M_2) \to \mathbb{P}(P_1)\to \mathbb{P}(P_0)
    \]
    is exact in $\modu_{\infty}-\mcP$. 
\end{itemize}
\end{lem}

\begin{proof}
(2) implies (1): If $M\in \modu_{\infty}-\mcP$, we choose $f\colon P_1\to P_0$ admitting a weak kernel such that $M\cong \coKer \mathbb{P}(f)$. Then choose the complex as in (2), to obtain that $M \in \Reso (\Bild \mathbb{P})$, and therefore $\Bild \mathbb{P}$ is finitely resolving. 
Conversely, let $f \colon P_0 \to P_1$ be a morphism as in (2). 
This implies $Z=\Ker \Hom_{\mcP} (-,f) \in \modu_{\infty}-\mcP$. Assume (1), i.e. $\Reso (\Bild \mathbb{P})= \modu_{\infty}-\mcP$. This means there exists an exact sequence 
\[ 0\to \mathbb{P} (M_n) \to \cdots \to \mathbb{P}(M_2) \to Z \to 0 \]
in $\modu_{\infty}-\mcP$ with $M_i \in \mcE$. Since $\mathbb{P}$ is fully faithful, the claim follows. 
\end{proof}

\section{The derived category of an exact category}

We recall from Buehler \cite{Bue}, section 9: Let $\mcA$ be an idempotent complete  additive category. 
One can define the category of chain complexes $Ch(\mcA)$ and the homotopy category $K(\mcA)$, whose objects equal those if $Ch(\mcA)$ but the morphisms are the quotient group given by chain maps modulo chain maps  homotopic to the zero map. $K(\mcA)$ has the structure of a triangulated category induced by strict triangles in $Ch(\mcA)$. A chain complex
\[
A=(\cdots \xrightarrow{d^{n-1}} A^n\xrightarrow{d^{n}} A^{n+1}\xrightarrow{d^{n+1}}\cdots ) \quad \in Ch(\mcA)
\]
is called \textbf{left bounded} if $A^n=0$ for $n<<0$, \textbf{right bounded} if $A^n=0$ for $n>>0$ and \textbf{bounded} if $A^n=0$ for $\lvert n \rvert >>0$. We denote by $K^+(\mcA), K^-(\mcA)$ and $K^b(\mcA)$ the full subcategories of $K(\mcA)$ whose objects are the left bounded, right bounded and bounded chain complexes respectively. By definition $K^b(\mcA)=K^+(\mcA)\cap K^-(\mcA)$. The subcategories $K^*(\mcA)$ with $*\in \{+,-,b\}$ are triangulated subcategories but not closed under isomorphism in $K(\mcA)$ unless $\mcA=0$.  \\

Now, let us assume that $\mcE=(\mcA, \mcS)$ is an exact category. 
A chain complex $A=(A^\bullet, d^\bullet)$ in $Ch(\mcA)$ is called \textbf{acyclic} if every differential factors as $d^n\colon A^n \to Z^{n+1}A\to A^{n+1}$ such that $Z^nA \to A^n \to Z^{n+1}A$ is an exact sequence (i.e. in $\mcS$).\\
Neeman proved that the mapping cone of a chain map between acyclic complexes is acyclic. This implies that the full subcategory $Ac(\mcE)=Ac_{\mcS}(\mcA)$ of $K(\mcA)$ given by acyclic complexes is a triangulated subcategory of $K(\mcA)$. If $\mcA$ is idempotent complete then $Ac(\mcE)$ is closed under isomorphism, every null-homotopic chain complex is acyclic and $Ac(\mcE)$ is a thick subcategory of $K(\mcA)$.
We define $Ac^*(\mcE):=Ac(\mcE)\cap K^*(\mcA)$ for $*\in \{+,-,b\}$. If $\mcA$ is weakly idempotent complete then the categories $Ac^*(\mcE)$ are thick subcategories of $K^*(\mcA)$ for $*\in \{+,-,b\}$. \\
Given any triangulated category $\mcC$ and thick subcategory $\mcT$ the \textbf{Verdier quotient} 
\[
\mcC/\mcT
\]
is defined via a localization (cf. \cite{Kr-book}, section 3.2, Prop. 3.2.2). It is again a triangulated category and the canonical functor $\mcC \to \mcC/\mcT$ is an exact functor. \\
As explained in \cite{Bue}, since $\mcA$ is idempotent complete, $Ac^*(\mcE)$ is a thick subcategory of $K^*(\mcA)$ for $*\in \{\pm, b, \emptyset\}$. 
The \textbf{derived category} of an exact category $\mcE=(\mcA, \mcS)$ is defined as the Verdier quotient 
\[
D(\mcE):= K(\mcA)/Ac (\mcE)
\]
Similarly, we define the \textbf{bounded/left bounded derived category} as the Verdier quotients 
\[
D^*(\mcE) := K^*(\mcA)/ Ac^*(\mcE), \quad *\in \{b, \pm \} 
\]
For more details and basic properties we refer to \cite{Kr-book}, chapter 4 and \cite{Bue}, section 10. 

Let $\mcE$ be an exact category with enough projectives $\mcP$. We define 
$K^{b,-}(\mcP)$ as the full subcategory of $K^-(\mcP)$ given by all complexes $X$ such that there exists an $n \geq 1$ such that $d^{n-1}$ and $d^{-n}$ are admissible and the truncated complexes \\
$\tau_{\leq (-n)}X =( \cdots \to X^{(-n)-1}\to X^{(-n)} \to \Ker d^{(-n)} \to 0 \cdots   )$, \\
$\tau_{\geq n}X=(\cdots 0 \to \coKer d^{n-1} \to X^{n+1} \to X^{n+2} \to \cdots )$\\ are acyclic, i.e. in $Ac(\mcE)\cap K^-(\mcP)$. Observe, that this depend on the ambient exact category $\mcE$, even though the notation suggests otherwise.
 
\begin{lem} (\cite{Kr-book}, Cor. 4.2.9)\label{noMorecrucialLemma}
Let $\mcE$ be an exact category with enough projectives $\mcP$. Then, there are triangle equivalences 
\begin{itemize}
    \item[(1)] $K^-(\mcP) \to D^-(\mcE)$
    \item[(2)] $K^{b,-} (\mcP) \to D^b(\mcE)$
\end{itemize}
\end{lem}

\begin{dfn}
Let $\mcE $ be an exact category and $\mcX$ a fully exact subcategory. We write $D^*(\mcX)$, $*\in \{b, +,-\}$ for its derived categories. There always exists an exact functor $D^*(\mcX) \to D^*(\mcE)$ induced by the inclusion $\mcX \to \mcE$ (because this is an exact functor). 
\end{dfn}



\begin{thm} \cite{HR}\label{HR-Thm}
If $\mcE$ is an exact category and $\mcX$ is a resolving or coresolving subcategory. 
\begin{itemize}
    \item[(i)] If $\mcX$ is resolving, the inclusion $\mcX \to \mcE$ induces a triangle equivalence \[D^-(\mcX) \to D^-(\mcE)\]
    if $\mcX$ is coresolving on $D^+()$. 
    \item[(ii)] If $\mcX$ is finitely resolving or finitely coresolving, the inclusion induces a triangle equivalence 
\[D^b(\mcX) \to D^b(\mcE).\] 
    \item[(iii)] If $\mcX$ is uniformley finitely resolving or uniformly finitely coresolving, it induces a triangle equivalence 
\[D(\mcX) \to D(\mcE).\] 
\end{itemize}
\end{thm}

\begin{rem} \label{ffres}
If $\mcE$ is an exact category with enough projectives and $\mcX$ a resolving subcategory, then the exact functor $\mcX\to \mcE$ induces a commuting diagram of exact functors
\[ 
\xymatrix{D^-(\mcX) \ar[r]&D^-(\mcE) \\
D^b(\mcX)\ar[r]\ar[u] & D^b(\mcE) \ar[u]}
\]
where the vertical arrows are the induced triangle functors by construction. By Lemma \ref{noMorecrucialLemma} they are fully faithful. By Theorem \ref{HR-Thm},(i), the upper horizontal functor is a triangle equivalence. This implies that the lower horizontal functor is fully faithful, too. 
\end{rem}

\section{Tilting theory in exact categories}

We define tilting subcategories and go through the list of results (cf. introduction) that we demand for tilting theory (i.e. induced derived equivalence, Aulsander correspondence, Bazzoni-like result and Brenner-Butler theorem). 

\subsection{Tilting subcategories}
\begin{dfn} \label{tiltingDef}
Let $\mcE=(\mcA, \mcS)$ be an exact category and $\mcT\subset \mcE$ a full subcategory. We call $\mcT$ a \textbf{tilting} subcategory if 
\begin{itemize}
    \item[(T1)] $\mcT$ is self-orthogonal and $\mcT^{\perp}$ has enough projectives which are given by $\mcP(\mcT^{\perp})=\mcT$\\
    (equivalently: $\mcT$ is closed under taking summands in $\mcE$ and $\mcT\subset \mcT^{\perp}\subset \gen (\mcT)$)
    \item[(T2)] $\mcT^{\perp}$ is finitely coresolving, i.e. it is coresolving and $\Cores (\mcT^{\perp})= \mcE$. 
\end{itemize}
By Prop. \ref{AB-thick}, if we assume (T1), we can replace (T2) by an equivalent:
\begin{itemize}
    \item[(T2')] $\Thick(\mcT^{\perp})=\mcE$.  
\end{itemize}

We call $\mcT$ $n$-\textbf{tilting} if it is tilting with $\Cores_n(\mcT^{\perp})=\mcE$. We call an object $T$ in $\mcE$ \textbf{tilting} if $\add(T)$ is a tilting subcategory. 
Assumption (T1) implies $\mcT$ and $\mcT^{\perp}$ are closed under summands and we have a well-defined fully faithful exact functor into infinitely presented $\mcT$-modules $f_{\mcT}\colon \mcT^{\perp} \to \modu_{\infty}-\mcT$, $X \mapsto \Hom_{\mcE}(-,X)\lvert_{\mcT}$ since $\mcT=\mcP(\mcT^{\perp})$ and $\mcT^{\perp}$ has enough projectives.\\  
Assumption (T2) implies that the inclusion $\mcT^{\perp}\subset \mcE$ induces a triangle equivalence $D^b(\mcT^{\perp}) \to D^b(\mcE)$.
Assumption (T1) and (T2) imply that the inclusion $\mcT \to \mcE$ gives rise to a triangle equivalence \[K^{b,-}(\mcT) \cong D^b(\mcT^{\perp}) \to D^b(\mcE). \]
\end{dfn}


Here is another way to express (T2) but we will need the assumption that the exact category has enough injectives. 

\begin{lem} \label{reformulate}
Let $\mcT$ be an additive  subcategory in an exact category $\mcE$. Then: 
    $\Cores (\mcT^{\perp})=\mcE$ (resp. $\Cores_n(\mcT^{\perp})=\mcE$) implies $\mcE= \bigcup_{n\geq 0} \mcT^{\perp_{\geq n+1}}$ (resp. $\pdim_{\mcE}\mcT\leq n$).\\ If $\mcE$ has enough injectives , then the converse is true. 
\end{lem}
\begin{proof} Observe that $A \in \Cores_n(\mcT^{\perp})$ implies $\Ext^{i}(T, A)=0, i\geq n+1$ for every $T\in \mcT$, i.e. $A\in \mcT^{\perp_{ \geq n+1}}$. If $\mcE$ has enough injectives and $A \in \mcT^{\perp_{\geq n+1}}$, then we have by dimension shift $\Omega^{-n}A \in \mcT^{\perp}$. So, the injective coresolution of $A$ truncated at the $n$-th cosyzygy shows that $A \in \Cores_n (\mcT^{\perp})$.
\end{proof}

\begin{cor}
Assume that $\mcT$ is a tilting subcategory of an exact category $\mcE$ and $n \geq 0$. Then the following are equivalent: 
\begin{itemize}
    \item[(1)] $\mcT$ is $n$-tilting
    \item[(2)] $\pdim_{\mcE} \mcT \leq n$
\end{itemize}
\end{cor}
\begin{proof}
By Lemma \ref{reformulate}, we have the implication (1) implies (2). Assume (2), we claim: $\Cores (\mcT^{\perp})\subset \Cores_n(\mcT^{\perp})$. Let $A\in \Cores (\mcT^{\perp})$, then there exists an exact sequence 
\[ 0\to A \to X_0 \to \cdots \to X_s \to 0, \quad X_i \in \mcT^{\perp}, 0\leq i \leq s  \] 
Wlog. assume that $s>n-1$. Let $Z_i=\Bild (X_i \to X_{i+1}), Z_{-1}:=A$.
Applying $\Hom(T,-)$, $T\in \mcT$, gives by dimension shift  $\Ext^j_{\mcE} (T,Z_{n-1})\cong \Ext^{j+1}_{\mcE}(T,Z_{n-2})\cong \cdots \cong \Ext^{j+n}_{\mcE} (T, A) =0$ for all $j\geq 1$ since $\pdim_{\mcE} \mcT \leq n$. So, $Z_{n-1} \in \mcT^{\perp}$. 
\end{proof}

\begin{lem} \label{weakBaz} Let $\mcT$ be an additive category of an exact category $\mcE$. 
 If $\mcT$ fulfills (T1) and $\pdim_{\mcE}\mcT\leq n$, then we have \[ \mcT^{\perp} =\gen_{n-1}(\mcT) = \pres_{n-1}(\mcT)=\gen_{\infty}(\mcT) =\pres_{\infty}(\mcT) .\]
\end{lem}

\begin{proof}
The assumption $\mcT^{\perp_{\geq n+1}}=\mcE$ together with $\mcT$ self-orthogonal implies by an easy dimension shift argument that $\pres_{n-1} (\mcT) \subset \mcT^{\perp}$. Now, by assumption (T1) we have $\mcT^{\perp}\subset \gen_{\infty}(\mcT)$ since $\mcT=\mcP(\mcT^{\perp})$ and $\mcT^{\perp}$ has enough projectives. Since $\gen_{\infty}(\mcT) \subset \pres_{n-1}(\mcT)$ trivially and $\gen_{n-1}(\mcT)$ and $\pres_{\infty}(\mcT)$ are intermediate between these two, the claim follows. 
\end{proof}

For infinitely generated tilting modules, the perpendicular category of a tilting module of finite projective dimension has been described in \cite{Baz}, Theorem 3.11. In \cite{Wei}, Theorem 1.1, this description has been proven for relative tilting modules over artin algebras. Therefore, we call the generalization to tilting subcategories in exact categories \emph{Bazzoni's description}. 

\begin{cor} \label{Baz}(\emph{Bazzoni's description})
Let $\mcT$ be an additive category of an exact category $\mcE$ which is closed under taking summands and assume $\pdim_{\mcE}\mcT\leq n$, then (T1) is equivalent to the following: 
  \[  \text{\emph{(T1')}}\quad \mcT\subset \mcT^{\perp}\subset \pres_{n}(\mcT)\]
and also to
\[  \text{\emph{(T1'')}}\quad \mcT\subset \mcT^{\perp}\subset \add(\pres_{n}(\mcT))\]
where for a full subcategory $\mcM$ of $\mcE$: $\add(\mcM)$ denotes the full subcategory of all summands of finite direct sums of objects in $\mcM$.\\ 
Furthermore, in this case we have $\pres_{n-1}(\mcT)=\pres_n(\mcT)=\mcT^{\perp}$.  
\end{cor}
\begin{proof}
The implication "(T1) implies (T1')" is given by Lemma \ref{weakBaz} using the inclusions $\pres_{\infty} (\mcT) \subset \pres_n(\mcT) \subset  \pres_{n-1}(\mcT)$. Also loc. cit. proves the statement $\pres_{n-1}(\mcT)=\pres_n(\mcT)$. The implication "(T1') implies (T1'')" is trivial since $\mcT^{\perp} = \add (\mcT^{\perp})$. \\
Assume (T1''), i.e. $\pdim_{\mcE}\mcT\leq n$ and $\mcT\subset \mcT^{\perp}\subset \add (\pres_{n}(\mcT) )$. To see (T1), it is enough to prove $\mcT^{\perp}\subset \gen (\mcT)$. Observe that $\pres_{n-1} (\mcT) \subset \mcT^{\perp}$ is fulfilled by Lemma \ref{basicProp}, (b). 
Let $X \in \mcT^{\perp}$, then there exists an $L \in \mcE$ such that $X\oplus L \in \pres_n(\mcT)$ - in particular $L \in \mcT^{\perp}$ and there is an exact sequence \[ 
0\to X_{n+1} \to T_{n} \xrightarrow{f_{n}} \cdots \xrightarrow{f_1} T_0 \to X \oplus L\to 0
\]
with $T_i \in \mcT$. 
Let $X_i:=\Bild (f_i), X\oplus L:=X_0$. For every $j\geq 1, T \in \mcT$ we have $\Ext^j (T, X_i) \cong \Ext^{j+1}(T, X_{i+1})$. This implies $X_1 \in \mcT^{\perp}$ since $\Ext^j(T,X_1) \cong \Ext^{n+j} (T, X_{n+1})=0$ by assumption. Therefore $T_0 \to X\oplus L$ is a deflation in $\mcT^{\perp}$. This implies that the composition with the projection onto the summand $T_0 \to X \oplus L \to X$ is a deflation in $\mcT^{\perp}$ - (since projections onto summands are deflation in every exact category and the composition of two deflations is a deflation).
Therefore, (T1) is fulfilled. 
\end{proof}

\begin{exa}
\begin{itemize}
    \item[(1)] If $\mcE$ is an exact category with enough projectives $\mcP$, then $\mcP$ is a tilting subcategory. In fact, this is the only $0$-tilting subcategory.
    \item[(2)] If $\mcT$ is a subcategory of an exact category $\mcE$ which fulfills (T1), then $\mcE'=\Cores (\mcT^{\perp})=\Thick (\mcT^{\perp})$ is a thick subcategory and $\mcT$ is a tilting subcategory of $\mcE'$.
    \item[(3)] A counterexample: Let $\mcT$ be a skelletally small additive category. 
    Let $\mcE=\modu_n-\mcT$ be the category of $n$-finitely presented $\mcT$-modules (i.e. for $n=0$ finitely generated, $n=1$ finitely presented, etc.). Assume that $\mcE$ does not have enough projectives. 
    $\mcE$ is a fully exact category of the abelian category $\Modu-\mcT$ which contains the finitely generated projectives $\mcP_{\mcT}:=\{\Hom_{\mcT}(-,T) \mid T \in \mcT\}$. We have $\mcP_{\mcT}^{\perp}=\mcE$ (so (T2) is fulfilled) but since $\mcE$ does not contain enough projectives, (T1) is not fulfilled and it is not a tilting subcategory. \\
    So, in our definition if you have a ring $R$ which is not left coherent, i.e. $\mcE=R-\modu_1$ is not abelian, then $R$ is not a tiliting module in $\mcE$. But $R$ is always a tilting module in $R-\modu_{\infty}$ by example (1). 
\end{itemize}
\end{exa}

\begin{lem}\label{thickTEqE}
Let $\mcE$ be an exact category and $\mcT$ an $n$-tilting subcategory. Then, all objects in $\mcT\subset \mcP^{<\infty}$ (i.e. all objects in $\mcT$ have finite projective dimension) and the following are equivalent: 
\begin{itemize}
    \item[(1)] $\mcE=\mcP^{<\infty}$ 
    \item[(2)] $\mcT^{\perp}\subset \mcP^{<\infty}$
    \item[(3)] $\mcT^{\perp} = \Reso (\mcT)$
    \item[(4)] $K^b(\mcT) =D^b(\mcT^{\perp}) $
    \item[(5)] $\mcE= \Thick (\mcT)$
    \item[(6)] $K^b(\mcT)=D^b(\mcE) $
\end{itemize}

\end{lem}
 \begin{proof}
 Since $\pdim_{\mcE}T\leq n$ for all $T \in \mcT$, $\mcT$ consists of ojects with finite projective dimension. This implies $\Thick (\mcT)\subset \mcP^{<\infty}$. Therefore (5) implies (1) trivially.
 Now assume (1). Then it follows $\mcT^{\perp}=\Reso (\mcT)=\Thick_{\mcT^{\perp}}(\mcT)$. This implies by (T2') that
 \[\mcE= \Thick (\mcT^{\perp})= \Thick (\Thick_{\mcT^{\perp}}(\mcT)) \subset \Thick (\mcT ). \]
 This means (1) and (5) are equivalent. The equivalences of (5) and (6) and the one of (3) and (4) follow from \cite{Kr-book}, Lemma 7.1.2 . The equivalence of (4) and (6) follows from the definition of a tilting subcategory. 
 The implications (1) implies (2), (2) implies (3) are clear. 
 \end{proof}

Out of curiosity, we also prove the following more general statement: 
\begin{lem} \label{thickT}
If $\mcE$ is an exact category and $\mcT$ is an $n$-tilting subcategory, then \[ \Thick (\mcT) = \mcP^{<\infty}\]
In particular, $\mcT$ is also an $n$-tilting subcategory of $\mcP^{<\infty}$. 
\end{lem}

\begin{proof}
Since $\mcT$ is $n$-tilting, we have $\mcT \subset \mcP^{<\infty}$. Since $\mcP^{<\infty}$ is a thick subcategory of $\mcE$, it follows that $\Thick (\mcT) \subset \mcP^{<\infty}$. \\
Now, let $X\in \mcP^{<\infty}$. There exists an exact sequence 
\[0\to X \to X_0 \to \cdot \to X_n \to 0\]
with $X_i \in \mcT^{\perp}$ by the axiom (T2). By Lemma \ref{weakBaz} we have $\mcT^{\perp}= \pres_{\infty}(\mcT)$, so we can apply Lemma \ref{technical-lemma} to find a short exact sequence $0 \to Z \to L \to X\to 0$ with $Z\in \mcT^{\perp}$ and $L \in \Cores (\mcT)$. 
This implies $L\in \Thick (\mcT)\subset \mcP^{<\infty}$ and $Z \in \mcP^{<\infty}$. 
Now clearly, $\mcT^{\perp}\cap \mcP^{<\infty} = \mcP^{<\infty}  (\mcT^{\perp}) \subset \Thick (\mcT )  $ since $\mcT^{\perp}$ has enough projectives given by $\mcT$. This implies $Z\in \Thick (\mcT)$ and therefore also $X\in \Thick (\mcT)$. \\
Since $\mcT^{\perp}\cap \mcP^{<\infty}$ has still enough projectives given by $\mcT$ we have (T1) trivially true. Secondly, since $\Thick_{\mcE} (\mcT)= \mcP^{<\infty}$, we have $\Thick_{\mcP^{<\infty}} (\mcT)= \mcP^{<\infty}$ and therefore (T2)' holds. 
\end{proof}
\begin{que}
Using the later theorem \ref{usualtilting} we can conclude that for exact categories with enough projectives $n$-tilting subcategories in $\mcE$ are precisely $n$-tilting subcategories in $\mcP^{<\infty}$. \\
Is this true for all exact categories? 
\end{que}

\subsubsection{The Auslander correspondence}
\label{Auslander-corresp}

Our definition has this trivial version of the Auslander correspondence as a consequence:
\begin{thm} (trivial Auslander correspondence) 
Let $\mcE$ be an exact category. The assignments $\mcT \mapsto \mcT^{\perp}$ and $\mcX\mapsto \mcP(\mcX)$ are inverse bijections between 
\begin{itemize}
    \item[(1)] the class of $\mcT$ (n-)tilting subcategory of $\mcE$
    \item[(2)] the class of $\mcX$ finitely coresolving subcategory (with $\Cores_n (\mcX) = \mcE$) which have enough projectives $\mcP=\mcP(\mcX)$ and which are of the form $\mcX= \mcP^{\perp}$
\end{itemize}
\end{thm}

Nevertheless, condition (2) can be reformulated 
in a more elegant way once we assume extra conditions on the exact category. This are the occurrences in the literature: 
\begin{itemize}
    \item[(*)] \cite{AR}, Theorem 5.5, original version for tilting modules in $\modu \La$ for an Artin algebra $\La$. 
    \item[(*)] \cite{Kr-book}, Theorem 7.2.18 in the case of a strongly homologically finite exact category $\mcE$ with enough projectives which are of the form $\mcP=\add (P)$ 
    for one object $P$. 
    \item[(*)] \cite{ZZ}, Theorem 4.15, for an extriangulated category with enough projectives and enough injectives. 
\end{itemize}
Here, we use that the different definitions of tilting are special cases of our definition (cf. section 5: Definitions of tilting). For example, the latest version of the Auslander correspondence is the following: 

\begin{thm} (Auslander correspondence) (\cite{ZZ}, Theorem 4.15)
Let $\mcE $ be an exact category with enough projectives and enough injectives. The assignments $\mcT \mapsto \mcT^{\perp}$ and $\mcX\mapsto {}^{\perp}\mcX\cap \mcX$ are inverse bijections between 
\begin{itemize}
    \item[(1)] the class of $\mcT$ (n-)tilting subcategory of $\mcE$
    \item[(2)] the class of $\mcX$ finitely coresolving subcategory (with $\Cores_n (\mcX) = \mcE$) which are covariantly finite and closed under summands.
\end{itemize}
\end{thm}

\begin{que} Can the previous result be proven only with the assumption that the exact category has enough projectives?
\end{que}


\section{Definitions of tilting}
\label{DefOfTilting}

In this section we look at our definition in exact categories with extra assumptions and compare this to existing definitions of tilting. Also, we discuss the relationship to tilting in triangulated subcategories in subsection \ref{tiltingInTri}.

\begin{pro} \label{ZZ-def}
If $\mcE$ is an exact category with enough injectives and $\mcT$ a full subcategory then $\mcT$ is $n$-tilting if and only if it fulfills the following two conditions:
\begin{itemize}
    \item[(i)] $\pdim_{\mcE}\mcT\leq n$
    \item[(ii)] $\mcT\subset \mcT^{\perp} \subset \pres_n (\mcT)$
\end{itemize}
Furthermore, in this case we have $\mcP(\mcE) \subset \Cores_n (\mcT)$.  
\end{pro}

\begin{proof}
Follows from Lem. \ref{Baz}, Lem. \ref{reformulate} and from Prop. \ref{AB-thick}. 
\end{proof}

\begin{rem}
This recovers the definition of $n$-tilting subcategories given in \cite{ZZ} (for exact categories with enough projectives and injectives). 
\end{rem}

\begin{thm} (and definition)\label{usualtilting}
Let $\mcE$ be an exact category with enough projectives and let $\mcT$ be a full additive subcategory. We define: 
\begin{itemize}
\item[(t1)] $\mcT$ is self-orthogonal and closed under summands. 
    \item[(t2)] $\pdim_{\mcE}\mcT \leq n$
    \item[(t3)] $\mcP(\mcE) \subset \Cores_n (\mcT)$
\end{itemize}
Then: $\mcT$ is $n$-tilting if and only if
$\mcT$ fulfills (t1), (t2), (t3).
\end{thm}

\begin{proof} 
Assume $\mcT$ is $n$-tilting. 
(T1) implies (t1), (t2) follows from Lemma \ref{reformulate}, (2), (a) and (t3) from Prop. \ref{AB-thick}, (2) with $\mcX=\mcT^{\perp}$.\\
Conversely, assume $\mcT$ fulfills (t1), (t2), (t3). It is straight forward to see that $\pres_n(\mcT) \subset \pres_{n-1} (\mcT)\subset \mcT^{\perp}$ by dimension shift using (t1) and (t2). \\
We claim: For every $X \in \mcE$ there is an exact sequence 
$X \to L \to Z$ in $\mcE$ such that $L \in \pres_n (\mcT)$ and $Z \in \Cores (\mcT)$. This short exact sequence implies $X \in \Cores (\mcT^{\perp})$. 
Using dimension shift and $\pdim_{\mcE} \mcT \leq n$ one can see $X \in \Cores_n (\mcT)$, so (T2) follows. \\
The claim follows using the beginning of a projective resolution of $X$
\[ 
P_n \to \cdots \to P_0 \to X \to 0
\]
with $P_i\in \mcP$. Recall by (t3) we have $P_i \in \Cores(\mcT)$, this implies that we can use Lemma \ref{technical} below to obtain the short exact sequence. \\
Now, assume $X \in \mcT^{\perp}$ and consider again the short exact sequence $X \to L \to Z$ with $L \in \pres_n(\mcT)$ and $Z \in \Cores (\mcT)$. Since $\mcT^{\perp}$ is inflation-closed we have $Z \in \mcT^{\perp}$. Now, it is easy to see that $\mcT^{\perp} \cap \Cores (\mcT) = \mcT$, so $Z \in \mcT$. Therefore we have $\Ext^1_{\mcE}(Z, X)=0$ which implies that the sequence splits
and we conclude that $X \in \add(\pres_n(\mcT))$. Thus we proved $\mcT^{\perp}\subset \add(\pres_n(\mcT))$. By Corollary \ref{Baz} we conclude that (T1) is fulfilled. 
\end{proof}

\begin{lem} \label{technical}
Let $\mcE$ be an exact category and $\mcT$ a full additive subcategory which is self-orthogonal and $n \geq 1$. Assume that we have an exact sequence 
\[
X_{n-1}\xrightarrow{f_{n-1}} X_{n-2} \to \cdots \to X_0 \xrightarrow{f_0} X\to 0
\]
with $X_i\in \Cores (\mcT )$. Then, there exists an exact sequence 
\[ T_{n-1} \xrightarrow{g_{n-1}} T_{n-2} \to \cdots \to T_0 \xrightarrow{g_0} L \to 0\]
with $T_i$ in $\mcT$ and an exact sequence $X\to L\to Z$ with $Z \in \Cores (\mcT)$. 
\end{lem}

The main ingredient to prove this is the following: The push out of an admissible morphism along an inflation is again admissible with the same kernel and cokernel as the admissible map that we pushed out (\cite{Bue}, Prop. 2.15)

\begin{proof}
 We choose a short exact sequence $X_{n-1}\to T_{n-1}\to Z_{n-1}$ 
with $T_{n-1} \in \mcT, Z_{n-1} \in \Cores (\mcT)$.
Then take the push out 
\[ 
\xymatrix{X_{n-1} \ar[r]^{f_{n-1}}\ar[d]& X_{n-2}\ar[d] \\
T_{n-1}\ar[r]^{h_{n-1}} & Y_{n-2}}
\]
since $\Cores(\mcT)$ is closed under extensions (by Lem \ref{basicProp} (a)) we have $Y_{n-2} \in \Cores(\mcT)$ since $\coKer (X_{n-2} \to Y_{n-2})=Z_{n-1}$. 
By construction $h_{n-1}$ admissible in $\mcE$ with kernel and cokernel of $h_{n-1}$ coincide with those of $f_{n-1}$. This means we constructed an exact sequence 
\[ 
T_{n-1} \to Y_{n-2} \to X_{n-3} \to \cdots \to X_0 \to X \to 0
\]
Now, we pick a short exact sequence $Y_{n-2} \to T_{n-2} \to Z_{n-2}$ with $T_{n-2} \in \mcT, Z_{n-2} \in \Cores (\mcT)$. Then we push out the admissible morphism $Y_{n-2} \to X_{n-3}$ along the inflation $Y_{n-2} \to T_{n-2}$ and proceed with the same method as before. 
Once $Y_0$ is constructed, choose the exact sequence $Y_0 \to T_0 \to Z_0$ with $T_0 \in \mcT, Z_0 \in \Cores (\mcT)$ and take the push out of the deflation $Y_0 \to X$ along the inflation $Y_0 \to T_0$. This gives another deflation $T_0 \to L$ such that $\Bild (T_1\to Y_0) = \Ker (T_0 \to L)$ and an inflation $X \to L$ with $\coKer (X \to L) = Z_0 \in \Cores (\mcT)$. 
\end{proof}

For an earlier lemma, we need the following modification of the previous lemma: 
Firstly, observe that we can replace $\Cores (\mcT)$ by $\copres_{\infty} (\mcT)$ as long as we know that $\copres_{\infty} (\mcT)$ is extension-closed. Secondly, if in the proof $f_{n-1}$ is an inflation, then we construct in the proof an exact sequence with $g_{n-1}$ an inflation. Thirdly, by passing to the opposite exact category, we have the following dual statement: 

\begin{lem} \label{technical-lemma}
Let $\mcE$ be an exact category and $\mcT$ a full additive subcategory which is self-orthogonal, $\pdim_{\mcE} \mcT <\infty$ and $n \geq 1$. Assume that we have an exact sequence 
\[
0\to X \to X_{0} \to \cdots \to X_{n-1}  \to 0
\]
with $X_i\in \pres_{\infty} (\mcT )$. Then, there exists an exact sequence 
\[ 0\to L \to T_0 \to \cdots \to T_{n-1} \to 0\]
with $T_i$ in $\mcT$ and an exact sequence $Z\to L\to X$ with $Z \in \pres_{\infty} (\mcT)$. 
\end{lem}

Observe that in the previous Lemma: The condition $\pdim_{\mcE} \mcT <\infty$ is needed to assure $\pres_{\infty} (\mcT)$ is extension-closed - cf. Lem. \ref{basicProp}.

\begin{rem}
In particular, this generalizes the usual definition of tilting modules of finite projective dimensions over an artin algebra (cf. Happel \cite{Hap}). Furthermore, it includes the generalization of Auslander and Solberg (cf. \cite{ASoII}) to so-called relative tilting modules (i.e. tilting objects in a different exact structure on the category of finitely presented modules over an artin algebra.  
 \end{rem}
 
Inclusion of perpendicular categories defines a partial order on all tilting subcategories on an exact category.
The previous theorem has the following application to this partial order. 
Here, given a fully exact subcategory $\mcX \subset \mcE$ we write $\Reso_{\mcX, m}(-)$ (resp. $\Cores_{\mcX, m}(-)$) for the category $\Reso_m(-)$ (resp. $\Cores_m(-)$) inside the fully exact category $\mcX$. 

\begin{pro} (and definition)
Let $\mcE$ be an exact category and $\mcT$ an $n$-tilting subcategory. Let $\tilde{\mcT}$ be another full subcategory which is self-orthogonal and closed under summands. 
Then the following are equivalent: 
\begin{itemize}
    \item[(a)] $\tilde{\mcT}$ is a $m$-tilting subcategory for some $m\geq 0$ of $\mcE $ and $\tilde{\mcT}^{\perp}\subset \mcT^{\perp}$.
    \item[(b)] $\tilde{\mcT}$ is a $m$-tilting subcategory for some $m\geq 0$ of $\mcE $ and $\tilde{\mcT}\subset \mcT^{\perp}$.
    \item[(c)]
    $\tilde{\mcT}$ is a $m$-tilting subcategory for some $m\geq 0$ of $\mcT^{\perp} $. 
    \item[(d)] There is an $m\geq 0$ such that  $\widetilde{\mcT}\subset \Reso_{\mcT^{\perp}, m}(\mcT) (\subset \mcT^{\perp})$ and $\mcT \subset \Cores_{\mcT^{\perp},m} (\tilde{\mcT})$.    
\end{itemize}
In this case, we write $\tilde{\mcT} \leq \mcT$.  
\end{pro}

Before, we give the proof. Let us note the following easy corollary: 
\begin{cor} \label{maxRigid}
Let $\mcE$ be an exact category and let $\mcT$ be an $n$-tilting subcategory and $\widetilde{\mcT}$ be an $m$-tilting subcategory for some $m\geq 0,n \geq 0$. If $\widetilde{\mcT} \subseteq \mcT$, then $\widetilde{\mcT}=\mcT$.
\end{cor}
\begin{proof}(of Cor. \ref{maxRigid}) 
The inclusion implies $\Ext^{>0} (\mcT, \widetilde{\mcT})=0=\Ext^{>0}(\widetilde{\mcT}, \mcT)$, therefore $\widetilde{\mcT}\leq \mcT \leq \widetilde{\mcT}$ and since $\leq$ is a partial order, they are equal. 
\end{proof}

\begin{proof}
We are going to show the equivalences of 1. (c) and (d), 2. (a) and (b) and then 3. (b) and (c). \\
1. The equivalence is a direct consequence of Theorem \ref{usualtilting} since $\mcT^{\perp}$ is an exact category with enough projectives. \\
2. Clearly (a) implies (b). So let us assume (b) and let $X \in \tilde{\mcT}^{\perp}$ and $T \in \mcT$. Since $\tilde{\mcT}^{\perp}$ has enough projectives given by $\tilde{\mcT}$, we find short exact sequences (in $\mcE$) setting $X:=X_0$
\[ 
0\to X_i \to \tilde{T}_i \to X_{i-1}\to 0
\]
with $\tilde{T}_i\in \tilde{\mcT}$ and $X_i \in \tilde{\mcT}^{\perp}$ for all $i \geq 1$. Then we have 
$\Ext^j_{\mcE} (T, X_{i-1}) \cong \Ext^{j+1} (T, X_i)$ for all $j\geq 1, i\geq 1$ implying $\Ext^j_{\mcE} (T, X_0) \cong \Ext^{j+n}_{\mcE} (T, X_n)=0$ since $\pdim_{\mcE} \mcT \leq n$. \\
3. Assume (b). Since (b) is equivalent to (a), we have  $\tilde{\mcT}^{\perp} \subset \mcT^{\perp}$. This implies that the perpendicular category of $\tilde{\mcT}$ inside $\mcT^{\perp}$ coincides with the perpendicular category $\tilde{\mcT}^{\perp}$ inside $\mcE$. This implies (T1). It is straight forward to see that $\Cores_{\mcT^{\perp}, m}(\tilde{\mcT}^\perp )= \Cores_m (\tilde{\mcT}^\perp) \cap \mcT^{\perp}$ holds and therefore (T2). \\
Assume (c). Since the equivalence to (d) implies $\mcT \subset \Cores_m(\tilde{\mcT})$ one can deduce $\tilde{\mcT}^{\perp} \subset \mcT^\perp$ (take an exact sequence $0\to T \to \tilde{T}_0 \to \cdots \to \tilde{T}_s \to 0$, $T \in \mcT, \tilde{T}_i \in \tilde{\mcT}$, take $X \in \tilde{\mcT}^\perp$ and apply $\Hom_{\mcE} (-, X)$ to the sequence to conclude $\Ext^{>0}_{\mcE}(T,X)=0$). This implies that the perpendicular category of $\tilde{\mcT}$ in $\mcT^\perp$ and in $\mcE$ coincide, therefore (T1) is fulfilled. To show (T2), observe that $\tilde{\mcT}\subset \Reso_m(\mcT)$ implies (since $\pdim_{\mcE} \mcT\leq n$) that there is a $t\geq 0$ such that $\pdim_{\mcE} (\tilde{\mcT}) \leq t$. Furthermore since $\mcT^{\perp}$ is coresolving and $\mcT^\perp \subset \Cores_m (\tilde{\mcT})$ we conclude that $\tilde{\mcT}^\perp$ is coresolving. This already implies that $\mcE = \Cores_t (\tilde{\mcT}^\perp )$, to see this: Let $X_0 \in \mcE$, then there are short exact sequences 
$0\to X_i \to Z_i \to X_{i+1}\to 0$ in $\mcE$ with $Z_i \in \tilde{\mcT}^\perp$, $i\geq 0$. Let $\tilde{T} \in \tilde{\mcT}$, then one has by dimension shift $\Ext^j_{\mcE} (\tilde{T}, X_t)\cong \Ext^{j+t}(\tilde{T}, X_0)=0$ for all $j \geq 1$ since
$\pdim\tilde{T} \leq t$. 
\end{proof}

Of course one has the dual result for cotilting subcategories (passing to opposite categories gives a poset isomorphism between $m$-tilting subcategories for some $m$ in $\mcE$ and $m$-cotilting subcategories for some $m$ in $\mcE^{op}$). 

\begin{rem}\label{mutation}
The previous proposition can be used to obtain the following (one-sided) \emph{mutation}: Given an $n$-tilting subcategory $\mcT=\mcM\vee \mcX$ with $\mcX \subset \gen(\mcM)$, define $\mcY= \Omega^-_{\mcM}\mcX$ and $\tilde{\mcT}=\mcM\vee \mcY$. If $\tilde{\mcT}$ is self-orthogonal (or equivalently: $\mcY \subset \cogen (\mcM), \mcX = \Omega_{\mcM} \mcY$), then $\tilde{\mcT}$ is $m$-tilting for some $m$ and $\tilde{\mcT}\leq \mcT$. \\
Just observe: $\tilde{\mcT} \subset \Reso_{\mcT^{\perp},1} (\mcT)\subset \mcT^{\perp}$ and $\mcT \subset \Cores_{\mcT^{\perp}, 1} (\tilde{\mcT})$. Therefore (d) in the previous Proposition applies. 
\end{rem}


Let $\mcE$ be an exact category and $n \geq 0$. 
Let $n-\tilt (\mcE)$ be the class of $n$-tilting subcategories.  

\begin{cor} (of Thm. \ref{usualtilting}) Let $\mcE$ be an exact category with enough projectives $\mcP$ and $\mcX$ a resolving subcategory of $\mcE$. Then one has 
\[n-\tilt (\mcX) = \{ \mcT \in n-\tilt (\mcE) \mid \mcT \subset \mcX \}.\]
In particular, we also have 
$n-\tilt ( \mcE) \cong \{ \mcT \subset n-\tilt (\modu_{\infty}-\mcP)\mid \mcT \subset \Bild \mathbb{P} \}, \; \mcT' \mapsto \mathbb{P} (\mcT') $
where $\mathbb{P} \colon \mcE \to \modu_{\infty} -\mcP, X \mapsto \Hom_{\mcE} (-,X)\lvert_{\mcP}$. 
\end{cor}

\begin{proof}
This is clear from the Theorem \ref{usualtilting}. The second statement follows since $\mathbb{P}$ is fully faithful, exact and $\Bild \mathbb{P}$ is a resolving subcategory of $\modu_{\infty}-\mcP$. 
\end{proof}

One of the obvious questions is when can one restrict to tilting objects and when is it necessary to study more general tilting subcategories? 
In short, at least if you have enough projectives and a Krull-Schmidt category then the category of projectives should tell you the answer. More detailed, we have: 
\begin{thm} \label{finProj} 
If $\mcE$ is a Krull-Schmidt, Hom-finite, exact category with enough projectives. 
\begin{itemize}
    \item[(1)] If there is a tilting object $T$, then we have $\lvert \mcP(\mcE) \rvert =\lvert T \rvert <\infty$.
    \item[(2)] If $\lvert \mcP (\mcE) \rvert < \infty$ and there is a tilting category $\mcT$, then we have $\lvert \mcT \rvert = \lvert \mcP(\mcE) \rvert < \infty$.
\end{itemize}
\end{thm}

\begin{proof} (of Thm \ref{finProj}) Let $\Gamma:= \End_{\mcE} (T)$. Since $\mcP=\mcP(\mcE)\subset \cogen^1_{\mcE} (T)$, we have by \cite{S}, Lemma 2.1 that $\Hom_{\mcE} (-,T) \colon \mcP \to \Gamma - \Modu$ is full and faithful. For every $P,P' \in \mcP$ we have $0 = \Ext^{i}_{\mcE} (P',P) \cong \Ext^i_{\Gamma} (\Hom_{\mcE}(P,T), \Hom_{\mcE} (P', T)) $ by \cite{S}, Lemma 3.3. Let $Q$ be the direct sum of all indecomposable projectives appearing in a minimal projective resolution of $T$. 
When we apply $\Hom_{\mcE} (-,T)$ to the projective resolution of $T$ and the exact sequence in (3), we conclude that $\Hom_{\mcE}(Q,T)$ is a tilting $\Gamma$-module. Since tilting modules are maximal rigid, we have $\add (\Hom_{\mcE} (Q, T))= \add ( \Hom_{\mcE} (\mcP, T))$. This implies $\lvert \mcP \rvert =\lvert Q \rvert= \lvert \Gamma \rvert = \lvert T \rvert$.
\end{proof}


 \subsection{Tilting in triangulated categories} \label{tiltingInTri}
 
 In this subsection we compare our definition of a tilting subcategory with the definition of a tilting subcategory in the bounded homotopy category of projectives and in the bounded derived category. 
 
 \begin{dfn}
 In a triangulated category $\mcC$ one defines $T\in \mcC$ to be tilting if $\Hom (T, \Sigma^i T)=0$ for $i\neq 0$ and the smallest thick subcategory that contains $T$ is $\mcC$. Similarly we call a full additive subcategory $\mcT$ of $\mcC$ \textbf{tilting} if
 \begin{itemize}
     \item[(Tr1)] $\Hom (T, \Sigma^i T')=0$ for $i\neq 0$ for all $T,T'$ in $\mcT$ and
     \item[(Tr2)] $\Thick_{\Delta} (\mcT)=\mcC$
     \end{itemize}
 \end{dfn}
 
 We recall the following Lemma which we want to use: 
 \begin{lem} (\cite{Kr-book}, Lemma 7.1.2) \label{Kr-lemma}
 Let $\mcX$ be a self-orthogonal subcategory in an exact category $\mcE$. We consider $\mcE \subset D^b(\mcE )$ as stalk complexes in degree zero.  
 Then the following are equivalent for an object $X\in \mcE$: 
 \begin{itemize}
     \item[(1)] $X\in \Thick (\mcX)$
     \item[(2)] $X\in \Thick_{\Delta} (\mcX)$
 \end{itemize}
 In particular, we have $\Thick (\mcX)=\mcE$ if and only if $\Thick_{\Delta} (\mcX) = D^b(\mcE)$. 
 \end{lem}

 \begin{lem} 
 If $\mcE$ is an exact category and $\mcT$ a full additive subcategory. Then the following are equivalent: 
 \begin{itemize}
     \item[(1)] $\mcT\subset \mcE \subset D^b(\mcE)$ already lies in $K^b(\mcP(\mcE))$ and gives rise to a tilting subcategory in $K^b(\mcP(\mcE))$. 
     \item[(2)] $\mcT$ is self-orthogonal and $\Thick (\mcP (\mcE))= \Thick (\mcT ) \subset \mcE$ 
     \item[(2')] $\mcT$ is self-orthogonal and 
     $\mcP(\mcE)\subset \Cores (\mcT)$, $\mcT \subset \Reso (\mcP(\mcE))$
 \end{itemize}
 \end{lem}
 Of course, in the situation that $\mcE$ has enough projectives and $\mcT=\add (T)$ contravariantly finite then (2') is equivalent to $\mcT$ being $n$-tilting (for some $n$).
 
 \begin{proof}
 We have for a self-orthogonal subcategory $\mcT \subset \mcE$: $\mcT \subset K^b(\mcP(\mcE))$ if and only if $\mcT \subset \Thick(\mcP(\mcE))=\Reso (\mcP(\mcE))$ by Lem. \ref{Kr-lemma} (\cite{Kr-book}, Lem 7.1.2). In this case, we have $\mcT$ self-orthogonal in $\mcE$ if and only if $\mcT$ fulfills (Tr1) in $K^b(\mcP(\mcE))$. Furthermore, we also have by loc.cit.:  
 $\Thick (\mcT)=\Thick (\mcP(\mcE))\subset \mcE$ if and only if
 $K^b(\mcP(\mcE))=\Thick_{\Delta }(\mcP(\mcE))=\Thick_{\Delta} (\mcT)\subset D^b(\mcE)$. \\
 We now observe that for a self-orthogonal category $\mcT$ we have $\mcP(\mcE)\subset \Thick (\mcT)$ implies $\mcP(\mcE)\subset \Cores (\mcT)$ by \cite{Kr-book}, Lem. 7.1.6. , this finishes the proof.
 \end{proof}

 \begin{lem} \label{Kr-def}
 If $\mcE$ is an exact category and $\mcT$ a full subcategory. Then, the following are equivalent
 \begin{itemize}
     \item[(1)] $\mcT$ is a tilting subcategory in $\mcE$ and $\Thick (\mcT)=\mcE $.
     \item[(2)] $\mcT$ is self-orthogonal and $\Thick (\mcT)=\mcE$
     \item[(2')] $\mcT$ is a tilting subcategory in $D^b(\mcE)$
 \end{itemize}
 \end{lem}
 
 \begin{proof} Clearly, (1) implies (2). 
 Assume (2), then $\mcT^{\perp}=\Reso (\mcT)$ by \cite{Kr-book}, Prop. 7.10. This implies that $\mcT= \mcP(\mcT^{\perp})$ and that $\mcT^{\perp}$ has enough projectives, so (T1). Also $\Thick (\mcT) = \mcE$ implies $\Cores (\mcT^{\perp}) = \Thick (\mcT^{\perp}) = \mcE$, so (T2) and $\mcT$ is tilting. The equivalence of (2) and (2') follows from Lem. \ref{Kr-lemma}. 
 \end{proof}
 
 \begin{rem}
 Lemma \ref{Kr-def}, (2) recovers the definition of tilting from \cite{Kr-book}, chapter 7 as a special case of our definition.  
 \end{rem}

\section{Induced triangle equivalences}
\label{Ideq}
\begin{dfn}
Let $\mcT$ be a tilting subcategory in an exact category $\mcE$. From now on, we call the functor 
 \[ 
 \boxed{f_{\mcT} \colon \mcT^{\perp} \to \modu_{\infty}-\mcT, \quad X\mapsto \Hom_{\mcE} (-,X)\vert_{\mcT}}
 \]
the \textbf{tilting functor} of $\mcT$. Furthermore, since $\mcT^{\perp}$ is finitely coresolving in $\mcE$ and $\Bild f_{\mcT}$ is resolving in $\modu_{\infty}-\mcT$, we have an exact functor 
\[ 
\boxed{F_{\mcT} \colon D^b(\mcE) \xrightarrow{\cong} D^b(\mcT^{\perp}) \xrightarrow{\cong} D^b(\Bild f_{\mcT} ) \to D^b(\modu_{\infty}-\mcT )}
\]
which we call the \textbf{derived tilting functor} of $\mcT$. 
\end{dfn}
Derived tilting functors are exact and fully faithful by remark \ref{ffres}

\begin{dfn}
Let $\mcT$ be a ($n$-)tilting subcategory in an exact category $\mcE$. We say that $\mcT$ is \textbf{ideq} ($n$-)-tilting if $F_{\mcT}$ is a triangle equivalence (i.e. essentially surjective).  We call $\mcT$ $m$-\textbf{ideq} ($n$-)tilting if $\Bild f_{\mcT}$ is $m$-resolving. 
\end{dfn}

\begin{rem} If $\Bild f_{\mcT}$ is finitely resolving in $\modu_{\infty}-\mcT$ then $F_{\mcT}$ is a triangle equivalence. 
By Lemma \ref{finResolving}, $\Bild f_{\mcT}$ is finitely resolving iff for every morphism $f \colon T_1 \to T_0$ in $\mcT$ which admits a sequence of successive weak kernels in $\mcT$ we find a complex in $\mcT^{\perp}$ 
    \[ 
    0\to M_m \to M_{m-1} \to \cdots \to M_2 \to T_1 \to T_0
    \]
    such that the induced sequence of restrictions of representable functors $0 \to \Hom (-,M_m)\lvert_{\mcT}\to \cdots \to \Hom (-,M_0)\lvert_{\mcT}\to \Hom_{\mcT}(-,T_1)\to \Hom_{\mcT}(-,T_0)$ in $\modu_{\infty}-\mcT$ is exact.  
    \end{rem}
\begin{rem}
If $\mcT$ is $m$-ideq $n$-tilting for some $n,m\geq 0$, then we have induced triangle equivalences also on the unbounded derived category
\[\mathbb{F}_{\mcT}\colon D(\mcE)\cong D(\mcT^{\perp}) \xrightarrow{f_{\mcT}} D(\Bild f_{\mcT}) \cong D(\modu_{\infty}-\mcT)
\]
\end{rem}


\begin{exa} 
Assume that $\mcT$ is a tilting subcategory such that every map $f\colon T_1 \to T_0$ in $\mcT$ admits a kernel in $\mcT^{\perp}$ (the monomorphism  on $\Ker f$ does not have to be an inflation), then $\mcT$ is $2$-ideq tilting
(just use the kernel of the map).\\
An instance of this is the following: Assume that $\mcT$ is $1$-tilting with $\mcT$ contravariantly finite in $\mcE$ (then: $\mcT^{\perp}= \gen (\mcT) = \pres (\mcT)$) and assume 
every morphism in $\mcT$ factors over a deflation in $\mcE$ (for example if $\mcE$ is abelian). Then every morphism in $\mcT$ has a kernel in $\mcT^{\perp}$. 
 \end{exa}

\subsection{Ideq tilting from equality with $\mcP^{<\infty}$}
Let us first observe the obvious: 
\begin{lem}
Let $\mcE$ be an exact category with enough projectives $\mcP$ and $m \geq 0$. Then, the following are equivalent: 
\begin{itemize}
    \item[(1)] $\mcE=\mcP^{<\infty}$ (resp. $\gldim \mcE\leq m< \infty$)
    \item[(2)] Every resolving subcategory is finitely resolving (resp. is $m$-resolving ). 
\end{itemize}
\end{lem}
\begin{proof}
(2) implies (1): $\mcP$ is a resolving subcategory, it is finitely resolving, i.e. $\Reso (\mcP) =\mcE$, (resp. $m$-resolving, i.e. $\Reso_m(\mcP)=\mcE$) if and only if $\mcP^{<\infty}=\mcE$ (resp. $\gldim \mcE\leq m<\infty)$. \\
(1) implies (2): If $\mcX$ is resolving then $\mcP \subset \mcX$. Therefore $\Reso (\mcP ) \subset \Reso (\mcX)$, $\Reso_m (\mcP )\subset \Reso_m(\mcX)$.
\end{proof}
Which brings us to this naive question:
\begin{que}
Is the previous lemma still true if we drop the assumption that $\mcE$ has enough projectives?
\end{que}

\begin{pro} \label{FinDerEq}
Let $\mcE$ be an exact category and $\mcT$ a tilting subcategory. 
If $\modu_{\infty}-\mcT = \mcP^{<\infty}$ (resp. $\gldim (\modu_{\infty}-\mcT)\leq m<\infty$), then $\mcT$ is ideq tilting (resp. $m$-ideq tilting) .  
\end{pro}

Another common triangle equivalence considered uses \emph{perfect complexes}:
\begin{dfn}
For an exact category $\mcE$ with enough projectives $\mcP$, 
we call $D^b_{perf}(\mcE):=K^b(\mcP)= \Thick_{\Delta}(\mcP)\subset D^b(\mcE)$ the triangulated subcategory of \textbf{perfect complexes}. 
\end{dfn}
Observe, that $D^b_{perf}(\mcE)=D^b(\mcE)$ is equivalent to $\mcE=\mcP^{<\infty}$ by Lemma \ref{thickTEqE} since $\mcP$ is a $0$-tilting subcategory. 

\begin{lem}\label{perfDer}
If $\mcE$ is an exact category with $\mcE=\mcP^{<\infty}$, and $\mcT$ tilting subcategory then we have an induced triangle equivalence  
\[ 
D^b(\mcE) \to D^b_{perf} (\modu_{\infty}-\mcT)
\]
\end{lem}

\begin{proof}
Since $\mcE=\mcP^{<\infty}$ implies $\mcT^{\perp}= \mcP^{<\infty}$ and since we have the additive equivalence $\mcP(\mcT^{\perp}) \to \mcP(\modu_{\infty}-\mcT), \; T \mapsto \Hom_{\mcT}(-,T)$ we obtain induced triangle equivalences 
\[ 
D^b(\mcE) \cong D^b(\mcT^\perp ) \cong K^b(\mcP(\mcT^{\perp})) \to K^b(\mcP(\modu_{\infty}-\mcT)) \cong D^b_{perf}(\modu_{\infty}-\mcT)
\]
\end{proof}

We also have the following: 
\begin{pro}\label{gldimIneq}
Let $\mcE$ be an exact category and $\mcT$ an $m$-ideq $n$-tilting subcategory. Then we have: 
\[\gldim \mcE \leq \gldim (\modu_{\infty}-\mcT) +n, \quad \gldim (\modu_{\infty}-\mcT) \leq \gldim \mcE +m  \]
\end{pro}

The proof follows directly from the following Lemma (and its dual statement). 
\begin{lem}\label{fingldim}
Let $\mcE$ be an exact category and $\mcX$ be a fully exact category. If we have $\Reso_n(\mcX)=\mcE$ for some $n \in \mathbb{N}$, then we have 
\[
\gldim \mcX \leq \gldim \mcE \leq \gldim \mcX +n
\]
\end{lem}

\begin{proof}
The first inequality is clear since $\Ext^i_{\mcX}=(\Ext^i_\mcE)\lvert_{\mcX}$. For the second, wlog. $\gldim \mcX=s\leq \infty$, let $E, L \in \mcE$, we claim $\Ext^{>(s+n)}_{\mcE}(L,E)=0$. By assumption, exists an exact sequence 
$0 \to X_n \to \cdots \to X_0 \to E\to 0$ with $X_i \in \mcX, 0 \leq i \leq n$. We apply $\Hom(X,-)$ with $X \in \mcX$ and obtain $\Ext^{>s}_{\mcE} (X,E)=0$. Now, we take the exact sequence $0\to Y_n \to \cdots \to Y_0 \to L \to 0$ with $Y_i \in \mcX, 0\leq i \leq n$, and apply $\Hom (-, E)$ and obtain $\Ext^{>(s+n)}_{\mcE}(L,E)=0$. 
\end{proof}

\subsection{Ideq tilting in exact categories with enough projectives}
The following is the most important result for this question (it is a generaliztion of Miyashita's theorem \cite{Mi}, Thm 1.16): 

\begin{thm} (Generalized Miyashita-Thm) \label{MiThm}
Let $\mcE$ be an exact category with enough projectives $\mcP$ and let $\mcT$ be an $n$-tilting subcategory which is essentially small. 
We consider the contravariant functor \[\Psi_{\mcT} \colon \mcE \to \mcT-\Modu ,\quad X\mapsto \Hom_{\mcE} (X,-)\lvert_{\mcT} \] and the covariant functor \[ \Phi_{\mcT} \colon \mcE \to \Modu -\mcT, \quad X\mapsto \Hom_{\mcE} (-,X)\lvert_{\mcT}.\] Let $\widetilde{\mcT}:= \Psi_{\mcT} (\mcP)$ and $\overline{\mcT}:=\Phi_{\mcP}(\mcT)$. Then we have: 
\begin{itemize}
    \item[(1)] $\widetilde{\mcT}$ is an $n$-tilting subcategory of $ \mcT-\modu_{\infty}$ and $\overline{\mcT}$ is an $n$-tilting subcategory of $\modu_{\infty}-\mcP$.
    \item[(2)] $\Psi_{\mcT}$ restricts to an equivalence $\mcP\cong \widetilde{\mcT}^{op}$ and $\Phi_{\mcP}$ to one $\mcT \cong \overline{\mcT}$. 
    \item[(3)] The category ${}_{\perp}\widetilde{\mcT} :=\{M \in \modu_{\infty}-\mcT\mid \Tor_{>0}^{\mcT} (M, \widetilde{\mcT} )=0 \}$ is a resolving subcategory of $\modu_{\infty}-\mcT$ with 
    $\Reso_n({}_{\perp}\widetilde{\mcT})= \modu_{\infty}-\mcT$.
    \item[(4)] The functor $\Phi = \Phi_{\overline{\mcT}} \colon \Modu-\mcP \to \Modu-\mcT, X\mapsto \Hom_{\Modu-\mcP}
    (\Phi_{\mcP}(-),X)\lvert_{\mcT}$ has a left adjoint $\Phi'\colon \Modu-\mcT \to \Modu-\mcP, X\mapsto (P \mapsto X\otimes_{\mcT} \Psi_{\mcT}(P))$. They restrict to inverse equivalences between
    \begin{itemize}
        \item[(i)] $\{M \in \Modu-\mcP \mid \Ext^{>0}_{\Modu-\mcP} (\mcT, M)=0\}$ and $\{ N \in \Modu- \mcT \mid \Tor_{>0}^{\mcT} (N, \widetilde{\mcT})=0\}$
        \item[(ii)] $\overline{\mcT}^{\perp} (\subset \modu_{\infty}-\mcP )$ and ${}_{\perp}\widetilde{\mcT} (\subset \modu_{\infty}-\mcT)$.
    \end{itemize}
    \item[(5)] We have a commutative triangle of exact functors (restricted to these subcategories)  
    \[ 
    \xymatrix{&\mcT^{\perp}\ar[dl]_{\Phi_{\mcP}}\ar[dr]^{\Phi_{\mcT}}& \\ \overline{\mcT}^{\perp}\ar[rr]_{\Phi} && {}_{\perp} \widetilde{\mcT} } 
    \]
\end{itemize}
In particular, $\overline{\mcT}$ is $n$-ideq $n$-tilting and we have an induced triangle equivalence $D^b(\modu_{\infty}-\mcP) \to D^b(\modu_{\infty}-\mcT)$. 
\end{thm}

Before we prove the previous theorem, let us state this Theorem as a corollary. 
\begin{thm}\label{ideqWep}
Let $\mcE$ be an exact category with enough projectives $\mcP$. Then the following are equivalent: 
\begin{itemize}
\item[(1)] $\mcE$ is equivalent as an exact category to a finitely resolving subcategory of $\modu_{\infty}-\mcP$
\item[(2)] There is an $n\in \mathbb{N}_0$ and an $n$-tilting subcategory of $\mcE$ which is ideq $n$-tilting.
\item[(3)] For every $n\geq 0$, every $n$-tilting subcategory of $\mcE$ is ideq $n$-tilting. 
\end{itemize}
\end{thm}

\begin{proof}
It is straight-forward to see that (1) is equivalent to $\mcP$ is ideq $0$-tilting. To prove the equivalences we show for a given $n$-tilting subcategory $\mcT$, we have: $\mcT$ is ideq $n$-tilting if and only if $\mcP$ is ideq $0$-tilting. For this it is enough to show (using Theorem \ref{ideqWep}) that we have a commutative diagram of triangle functors 
\[
\xymatrix{ & D^b(\mcE)\ar[dr]\ar[dl] & \\
D^b(\modu_{\infty}-\mcP) \ar[rr]&& D^b(\modu_{\infty}-\mcT)}
\]
where $D^b(\mcE)\to D^b(\modu_{\infty}-\mcP)$ is the derived functor of $\Phi_{\mcP}$ and $D^b(\mcE)\cong D^b(\mcT^{\perp}) \to D^b(\modu_{\infty}-\mcT)$ is the derived functor of $\Phi_{\mcT}|_{\mcT^{\perp}}$ and $D^b(\modu_{\infty}-\mcP) \to D^b(\modu_{\infty}-\mcT)$ is the triangle equivalence induced by the equivalence $\overline{\mcT}^{\perp} \to {}_{\perp}\widetilde{\mcT}$. But this follows immediately from loc. cit. (5). 
\end{proof}

We can prove the even stronger corollary of Theorem \ref{ideqWep}.

\begin{dfn}
We define $\Tilt (\mcE):=\bigoplus_{n\geq 0} n-\tilt (\mcE)$. The relation $\leq $ from Lemma .. defines a poset structure on this set.\\ 
Let $\mcE$ be an exact category. We say $\mcE$ is \textbf{tilting connected} if in the poset $\Tilt (\mcE)$ is non-empty and for every two element $\mcT$ and $\mcT'$ there is a finite sequence $\mcT_0=\mcT, \mcT_1 , \ldots , \mcT_r=\mcT'$ with $\mcT_i\leq \mcT_{i+1}$ or $\mcT_i \geq \mcT_{i+1}$, $0\leq i \leq r-1$. 
\end{dfn}

\begin{exa}
If $\mcE$ has enough projectives $\mcP$ then $\mcE$ is tilting connected since $\mcP$ is a globales maximum. \\
If the injectives $\mcI$ in $\mcE$ happen to be an $n$-tilting subcategory then $\mcE$ is tilting connected since $\mcI$ has to be a global minimum.  
\end{exa}

\begin{que}
Are exact categories are always tilting connected?
\end{que}

\begin{cor}
Let $\mcE$ be an exact category and $\mcT$ an $n$-tilting subcategory. Then the following are equivalent: 
\begin{itemize}
    \item[(1)] $\mcT$ is ideq $n$-tilting
    \item[(2)] Every $m$-tilting subcategory $\mcL$ with $\mcL\leq \mcT$ or $\mcT\leq \mcL$ (i.e. $\mcL$  comparable to $\mcT$) is ideq $m$-tilting. 
\end{itemize}
In particular, every connected component in $\Tilt (\mcE)$ is either ideq tilting (i.e. every $n$-tilting subcategory in it is ideq $n$-tilting) or it is not ideq tilting. 
\end{cor}

The proof of the main result uses an auxiliary preprint of the author \cite{S} in which some technical results are explained.  
\begin{proof} (of Theorem \ref{MiThm})
\begin{itemize}
    \item[(1)] We want to see that $\widetilde{\mcT}$ satisfies (t1),(t2) and (t3): Since $\mcT$ is an $n$-tilting subcategory of $\mcE$ we have by (t2) for every $T \in \mcT$ an exact sequence 
    \[ 
    0\to P_n \to \cdots \to P_0 \to T \to 0
    \]
    with $P_i \in \mcP$. We apply $\Psi_{\mcT}$ to it and obtain a complex 
    \[ 
    0\to \Psi_{\mcT} (T) \to \Psi_{\mcT} (P_0) \to \cdots \to \Psi_{\mcT} (P_n) \to 0 
    \]
    Since $\mcT$ fulfills (t1), it follows that this is exact, so $\widetilde{\mcT}$ fulfills (t3).\\
    Now, by (t3) for $\mcT$ we have for every $P \in \mcP$ an exact sequence \[0\to P \to T_0 \to \cdots \to T_n \to 0, \]
    we apply again $\Psi_{\mcT}$ and again by (t1) for $\mcT$ we get an exact sequence \[ 
    0\to \Psi_{\mcT} (T_n) \to \cdots \to \Psi_{\mcT} (T_0) \to \Psi_{\mcT}(P) \to 0
    \]
    which shows that $\widetilde{\mcT}$ fulfills (t2). Furthermore, $\mcP \subset \cogen^{\infty}_{\mcE} (\mcT)$ implies by \cite{S}, Lemma 3.3. that $\Ext^i(\Psi_{\mcT}(P), \Psi_{\mcT}(P'))=0$ for $P,P'\in \mcP$ and $0<i <\infty $, so (t1) holds for $\widetilde{\mcT}$.\\
    \noindent
    The second claim follows since $\Phi_{\mcP}$ is exact, fully faithful and preserves extension groups and maps projectives to projectives. 
    \item[(2)] Since $\mcP \subset \cogen^1_{\mcE} (\mcT)$ it follows that $\Psi_{\mcT}$ restricted to $\mcP$ is fully faithful by \cite{S}, Lemma 2.1. Since $\Phi_\mcP$ is fully faithful, the second claim is clear. 
    \item[(3)] By the properties of $\Tor$ the category   ${}_{\perp}\widetilde{\mcT}$ contains the projectives, is extension closed and deflation-closed, so it is resolving. \\
    The last statement follows when we consider the first $n$ terms of a projective resolution of $X\in \modu_{\infty}-\mcT$ 
    \[ 
    0\to \Omega^n \to T_{n-1} \to \cdots \to T_0 \to X \to 0
    \]
    with $T_i \in \mcP(\modu_{\infty}-\mcT)$. We claim that $\Tor_{>0}^{\mcT} (\Omega^n, \widetilde{\mcT} )=0$. By dimension shift $\Tor_{i}^{\mcT} (\Omega^n, \widetilde{\mcT} )=\Tor_{i+n}^{\mcT} (X, \widetilde{\mcT} )=0$ since $\pdim \widetilde{\mcT} \leq n$. 
    \item[(4)] It is standard to see that these functors form an adjoint pair (it should be seen as a Hom-Tensor adjunction), cf. \cite{S}, Lemma 3.7. \\
    (i) This is a straight forward generalization of the orginal result \cite{Mi}, Thm 1.16. We just mention it for completeness. \\
    (ii) 
    We want to see that both functors restrict to functors as claimed and that they are both fully faithful. \\
    By Lemma \ref{Baz} and \cite{S}, Lem. 3.13, Rem. 3.14 we have 
   \[ 
   \begin{aligned}
   \overline{\mcT}^\perp &=\gen_{\infty}^{\modu_{\infty}-\mcP}(\overline{\mcT})= \pres_{\infty}^{\modu_{\infty}-\mcP}(\overline{\mcT})\\
   &=\{ X \in \modu_{\infty}-\mcP \mid \varphi_X \text{ isom}, \Phi (X) \in \modu_\infty-\mcT, \Tor_{>0}^\mcT (\Phi(X), \widetilde{\mcT})=0\} 
   \end{aligned}
   \]
    therefore, the functor $\Phi$ restricts as claimed and is fully faithful since $\gen_{\infty}(\overline{\mcT})\subset \gen_1(\overline{\mcT})$ (again using \cite{S}, Lemma 1.1.) \\
    We are going to proof the following claim: 
    \begin{itemize}
        \item[(*)] 
        $\pres_{\infty}^{\modu_{\infty}-\mcP} (\overline{\mcT})=\pres_{\infty}^{\Modu-\mcP}(\overline{\mcT})$
        \end{itemize}
    
    \begin{itemize} 
    \item[proof of (*):]
    Given an exact sequence $\cdots \to T_m \to \cdots \to T_0 \to X \to 0$ in $\Modu-\mcP$ with $T_i \in \overline{\mcT}$, we claim $X \in \modu_{\infty}-\mcP$. 
    Now, we see this as a quasi-isomorphism of complexes (with terms in $\Modu-\mcP)$
     \[ \xymatrix{
     T_* \colon\quad  \cdots \ar[r] &T_m\ar[d] \ar[r] & \cdots \ar[r] &T_1\ar[r]\ar[d] &T_0\ar[r]\ar[d] & 0 \\
     X\colon \quad \cdots \ar[r] & 0 \ar[r] & \cdots \ar[r]& 0\ar[r]  & X\ar[r] & 0 } 
    \]
    Since $T_*$ is a complex in $\modu_{\infty}-\mcP$, we find by \cite{Bue}, Thm 12.7, a quasi-isomorphism $P_* \to T_*$ with $P_*$ is a complex $\cdots \to P_m \to P_{m-1}\to \cdots \to P_0 \to 0 $, $P_{-n}=0$ for all $n> 0$ with terms in $\mcP$, here a quasi-isomorphism means that the mapping cone of $P_* \to T_*$ is acyclic. Since composition of quasi-isomorphisms are quasi-isomorphsims, we have a quasi-isomorphism $P_* \to  X$ which means that the mapping cone which is a complex $ \cdots \to P_1 \to P_0 \to X \to 0$ is exact and therefore $X \in \modu_{\infty}-\mcP$. 
    \end{itemize} 
    Let $Y\in {}_{\perp}\widetilde{\mcT}$.  
    Since $Y\in \modu_{\infty}-\mcT$, e.g. there exists an exact sequence 
    \[ 
    \cdots \to \Phi(T_m) \to \cdots \to \Phi(T_0) \to Y \to 0
    \]
  with $\Phi(T_i) \in \mcP(\modu_{\infty}-\mcT)= \Phi (\overline{\mcT})$. Applying $\Phi'$ and using that $\Phi'\Phi (T) \cong T$ for all $T \in \overline{\mcT}$ yields a complex 
  \[ 
  \cdots \to T_m \to \cdots \to T_0 \to \Phi' (Y) \to 0
  \]
with $T_i \in \overline{\mcT}$. Since $\Tor_{>0}^{\mcT}(Y, \widetilde{\mcT})=0$ and $\Phi'$ is right exact, this complex is exact. This implies by (a) that $\Phi'(Y) \in \overline{\mcT}^{\perp}=\gen_{\infty}^{\modu_{\infty}-\mcP}(\overline{\mcT})$. Therefore, $\Phi'$ restricts as claimed. To see that it is fully faithful, we also have the consequence that applying $\Phi$ is again exact on this complex
\[ 
\cdots \to \Phi(T_m) \to \cdots \to \Phi(T_0) \to \Phi (\Phi' (Y)) \to 0
\] 
By the triangle identity of the adjunction, this implies that the unit $Y \to\Phi \Phi'(Y)$ is an isomorphism.\\
Now, an adjunction with unit and counit isomorphisms (i.e. fully faithful left and right adjoint) is an equivalence of categories. \item[(5)] It follows from the definition that $\Phi\circ \Phi_{\mcP}= \Phi_{\mcT}$. Since $\Phi_{\mcP}$ preserves all extension groups, it is clear that $\Phi_{\mcP}$ restricts to a functor as claimed, for $\Phi$ this has been proven in (4). 
\end{itemize}
\end{proof}

\begin{exa}
If $\mcE$ is an abelian category with enough projectives, then every $n$-tilting subcategory is ideq $n$-tilting. For example, let $\mcP$ be essentially small, then $\mcE=\Modu-\mcP$ is an abelian category with enough projectives given by summands of arbitrary direct sums of $\mcP$, i.e. $\mcP(\Modu-\mcP)=\rm{Add}(\mcP)$. It implies that the functor $\mathbb{P} \colon \mcE\to \modu_{\infty}-(\Add(\mcP))$  is an exact equivalence.
\end{exa}

\begin{exa}
Let $\mcE:= \mcP(\modu_{\infty}-\mcP )\subset \modu_{\infty}-\mcP$. We consider $\mcE$ as a fully exact subcategory of $\modu_{\infty}-\mcP$, then it is a resolving subcategory. Since the category $\mcE$ is semi-simple, it has a unique tilting subcategory $\mcT=\mcP(\mcE)=\mcE$ which is a $0$-tilting subcategory. Now, it follows from Thm. \ref{ideqWep} that $\gldim (\modu_{\infty}-\mcP )=m<\infty$ (resp. $\modu_{\infty}-\mcP= \mcP^{<\infty}$) if and only if $\Reso_m(\mcE)=\modu_{\infty}-\mcP$ (resp. $\Reso (\mcE ) = \modu_{\infty}-\mcP$) if and only if $\mcE$ is $m$-ideq $0$-tilting subcategory (resp. ideq $0$-tilting) of itself. 
\end{exa}
\subsection{Conjectural generalizations to ideq tilting in arbitrary exact categories.}
Here, we come accross the following unresolved question: 

\textbf{Open question:} Is there an analogue of Rickard's Morita theory for functor categories as follows: 

\begin{conj} (Strong Rickard conjecture) Let $\mcX, \mcY$ be two essentially small idempotent complete additive categories. Then the following equivalent:\\
\begin{itemize}
\item[(0)] $D(\Modu -\mcX) $ and $D(\Modu -\mcY)$ are triangle equivalent. 
\item[(1)] $D^-(\Modu-\mcX)$ and $D^-(\Modu-\mcY)$ are triangle equivalent.
    \item[(2)] $D^-(\modu_{\infty}-\mcX)$ and $D^-(\modu_{\infty}-\mcY)$ are triangle equivalent.
    \item[(3)] $D^b(\modu_{\infty}-\mcX)$ and $D^b(\modu_{\infty}-\mcY)$ are triangle equivalent.
    \item[(4)] $K^b(\mcX)$ and $K^b(\mcY)$ are triangle equivalent.
    \item[(5)] There is a tilting subcategory $\mcT$ of $K^b(\mcX)$ such that $\mcT\cong \mcY$ as additive categories. 
\end{itemize}
If we do not assume $\mcX$ and $\mcY$ to be small, then (2),(3),(4), (5) are still equivalent.
Furthermore, every triangle equivalence in (0)-(3) restricts to a triangle equivalence as in (4). 
\end{conj}

For $\mcX$ and $\mcY$ of finite type and idempotent complete, part of the conjecture is  
Rickard's Morita theorem for derived categories of rings (cf. \cite{Ric}, Thm 6.4, Prop. 8.1 - for example (3) implies (4) is proven in loc. cit only for right coherent rings - yet (4) implies (3) is proven for arbitrary rings). An alternative proof can be found in \cite{Kr-book}. For more general small additive categories partial answers are given by Keller \cite{K-DG}, Corollary in 9.2 and  Asadollahi-Hafezi-Vaheed \cite{AHV}, Thm 3.21.  \\
For our purpose, we need the following weaker statement:

\begin{conj} \label{RicLem} (Rickard-Lemma) Let $\mcX$, $\mcY$ be small idempotent complete additive categories. If $K^b(\mcX)$ and $K^b(\mcY)$ are triangle equivalent, then there exists  
a triangle equivalence $D^b(\modu_{\infty}-\mcX)\to D^b(\modu_{\infty}-\mcY)$ which restricts to a triangle equivalence $K^b(\mcX) \to K^b(\mcY)$.
\end{conj}

For $\add \La$ with $\La$ an arbitrary ring, the existence is proven in \cite{Ric} after Prop 8.1. One proof of the Rickard-Lemma should be as follows: \\
Let $\mcT\subset K^b(\mcY)$ be the image of $\mcX$ under the assumed triangle equivalence. We define $\Add(\mcY):=\mcP (\Modu-\mcY)$.
Then one shows that $\mcT\subset K(\Add(\mcY))$ fulfills (P1),(P2),(P3) in \cite {AHV} and the acyclic complexes in $K(\Add (\mcY))$ coincide with the $\mcT$-acyclic complexes in loc. cit (since $\Thick_{\Delta}(\mcT)$ in $K^b(\mcY)$ equals $K^b(\mcY)$)). Therefore, \cite{AHV}, Thm 3.21 can be applied to obtain a triangle equivalence $K^-(\Add(\mcX))\cong K^-(\Add(\mcT)) \to K^-(\Add(\mcY))$ which restricts to a triangle equivalence $K^b(\mcX) \to K^b(\mcY)$. Then, arguments of \cite{Ric} should generalize to intrinsic characterizations of the subcategory inclusions $K^{b,-}(\mcX) \subset K^-(\mcX) \subset K^-(\Add (\mcX))$ - these imply the claimed restricted triangle equivalence $D^b(\modu_{\infty}-\mcX) \to D^b(\modu_{\infty}-\mcY)$. Unfortunately, these results are not easy to puzzle together, so we leave this as a conjecture. 

\begin{cor} (of Rickard-Lemma) Assume that the Conjecture \ref{RicLem} holds. 
Let $\mcE$ be an exact category. Then for every  $n$-tilting and $m$-tilting subcategories $\mcT$ and $\mcT'$ which are small there exists a triangle equivalence $D^b(\modu_{\infty}-\mcT) \to D^b(\modu_{\infty}-\mcT')$
\end{cor}

\begin{proof}
We have the two derived tilting functors 
\[
\xymatrix{
D^b(\modu_{\infty}-\mcT) & D^b(\mcE) \ar[r] \ar[l] & D^b(\modu_{\infty}-\mcT')}
\]
Now, the thick subcategory of $D^b(\mcE)$ that $\mcT$ and $\mcT'$ generate is equal to 
\[ 
\Thick_\Delta (\mcT) = \Thick_{\Delta} (\mcP^{<\infty}) =\Thick_{\Delta} (\mcT')
\]
by Lemma \ref{thickT}. This implies that derived tilting functors restrict to triangle equivalences
\[
\xymatrix{
K^b(\mcT) &  \Thick_{\Delta} (\mcP^{<\infty})\ar[r] \ar[l] & K^b(\mcT')}
\]
Then the claim follows from the conjecture \ref{RicLem}. 
\end{proof}

\begin{rem}
If we could prove the stronger statement that every triangle equivalence $K^b(\mcT)\to K^b(\mcT')$ can be extended to a triangle equivalence as in the previous corollary, then we would obtain that the existence of one ideq $n$-tilting subcategory is equivalent to that all $m$-tilting are ideq $m$-tilting for all $m \geq 0$. This extension property for arbitrary triangle equivalences $K^b(\mcX)\to K^b(\mcY)$ would imply the following conjecture.
\end{rem}
\begin{conj} (Generalization of Thm. \ref{ideqWep}) Let $\mcE$ be an exact category and assume that there exists an $n \geq 0$ such that there is at least one $n$-tilting subcategory. Then the following are equivalent: 
\begin{itemize}
    \item[(a)] There is an $n\geq 0$ and an $n$-tilting subcategory which is ideq $n$-tilting
    \item[(b)] For every $m\geq 0$  every $m$-tilting subcategory is ideq $m$-tilting. 
    \item[(c)] There is a triangle equivalence $D^b(\mcE) \to D^b(\modu_{\infty}-\mcS)$ for some idempotent complete additive category $\mcS$ which restricts to a triangle equivalence $\Thick_{\Delta}(\mcP^{<\infty}) \to K^b(\mcS)$. 
\end{itemize}

\end{conj}
\begin{rem}
If every triangle equivalence $D^b(\mcE) \to D^b(\modu_{\infty}-\mcS)$ for some idempotent complete additive category $\mcS$ restricts to a triangle equivalence $\Thick_{\Delta}(\mcP^{<\infty}) \to K^b(\mcS)$, then the previous conjecture says for an exact category with at least one $n$-tilting subcategory: 
$n$-tilting $=$ ideq $n$-tilting for every $n \geq 0$ is equivalent to $\mcE$ is bounded derived equivalent to a category $\modu_{\infty}-\mcS$.  
\end{rem}

\subsection{When is the image of the tilting functor the perpendicular of a cotilting subcategory?}

Enomoto characterized in \cite{E-master}, Thm 2.4.11 when an exact category $\mcE$ is equivalent to the perpendicular category ${}^{\perp}\mcC$ of an $m$-cotilting subcategory $\mcC$ inside a functor category $\modu_{\infty}-\mcP$ (Clearly a perpendicular category of a cotilting subcategory is a finitely resolving subcategory). Observe that a necessary condition is that such an $\mcE$ has enough projectives and enough injectives. \\
Here Enomoto's notion of \emph{higher kernels} in an additive category plays a crucial role. 
\begin{dfn} (cf. \cite{E-master}, Def. 2.4.5) Let $\mcC$ be an additive category and $n \geq 1$, then we say that $\mcC$ has $n$-\textbf{kernels} if for every $f\colon C_1\to C_0$ in $\mcC$ there is a complex $0\to C_{n+1} \to \cdots \to C_2 \to C_1 \xrightarrow{f} C_0$ in $\mcC$ such that 
\[ 0 \to \Hom_{\mcC} (-, C_{n+1}) \to \cdots \to \Hom_{\mcC}(-,C_0)\]
is exact in $\modu_{\infty}-\mcC$. \\
If $\mcC$ is additionally an exact category, we say that $\mcC$ has $0$\textbf{-kernels} if every morphism $f$ in $\mcC$ can be factored as $f=id$ for a deflation $d$ and a monomorphism $i$. We say $\mcC$ has  $(-1)$\textbf{-kernels} if it is abelian. 
\end{dfn}

\begin{exa}
If $\mcT$ is ($n$-)tilting and has $m$-kernels for some 
$m \geq 1$, then $\mcT$ is $(m-1)$-ideq ($n$-)tilting (check the definitions).  
\end{exa}

\begin{pro} (\cite{E-master}, Pro. 2.4.6) 
Let $\mcC$ be essentially small additive category with weak kernels and $n \geq 1$. Then: $\modu_{\infty}-\mcC=\mcP^{\leq n+1}$ if and only if $\mcC$ has $n$-kernels. 
\end{pro}


\begin{thm}(Enomoto's Theorem, \cite{E-master}, Thm 2.4.11)
Let $\mcE$ be an idempotent complete, exact category and $\mcT$ a tilting subcategory and $m\geq 0$.
We write $f_{\mcT} \colon \mcT^{\perp} \to \modu_{\infty}-\mcT$ for the tilting functor $X\mapsto \Hom_{\mcE}(-,X)|_{\mcT}$. 
The following are equivalent:
\begin{itemize}
    \item[(1)] $\mcT$ has weak kernels and is $(m-1)$-ideq tilting and $\mcT^{\perp}$ has enough injectives   
    \item[(2)] There is an $m$-cotilting subcategory $\mcC$ in $\modu_{\infty}-\mcT$ with ${}^{\perp}\mcC = \Bild f_{\mcT}$
    \item[(3)] $\mcT^{\perp}$ has enough injectives and $(m-1)$-kernels
    \item[(4)] $\mcT^{\perp}$ has enough injectives and there is a category $\mcT \subset \mcM \subset \mcT^{\perp}$ which has $(m-1)$-kernels
    \end{itemize}
     In this case, $\Bild f_{\mcT}={}^\perp \mcC$ with $\mcC=f_{\mcT}(\mcI(\mcT^{\perp}))$ is $m$-cotilting.    
\end{thm}    

All examples that I know of this situation follow from Auslander's notion of a dualizing $R$-variety. 
\begin{dfn} (and Lemma)
Let $R$ be a commutative ring such that there is a duality $\kdual $ on finite length $R$-modules (e.g. if $R$ is a field).
Let $\mcA$ be an additive $R$-category such that $\Hom_{\mcA} (X,Y)$ is a finite length $R$-module for all $X,Y$ in $\mcA$. Then, $\kdual \colon \ModmcA \to \mcA-\Modu, F \mapsto F \circ \kdual$ is a duality (if $\mcA$ is essentially small). 
$\mcA$ is called a \textbf{dualizing R-variety} if the $F \mapsto F\circ D$ defines a duality between  finitely presented left and right $\mcA$-modules, i.e. a contravariant equivalence  
$\kdual \colon \modu_1-\mcA \to \mcA-\modu_1 \colon \kdual$. \\
In this case, 
$\modu_1-\mcA=\modu_{\infty}-\mcA$ is an abelian category with enough injectives and projectives and $\mcA^{op}$ is also a dualizing R-variety. 
\end{dfn}
    

\begin{cor} (of Miyashita's and Enomoto's Theorem) 
Let $\mcE$ be an exact category with enough projectives $\mcP$. Let $\mcT$ be an $n$-tilting subcategory of $\mcE$ and assume that there is a duality (i.e. contravariant equivalence) $\kdual \colon \modu_{\infty}-\mcT \to \mcT-\modu_\infty\colon \kdual $. Let $\widetilde{\mcT}=\Psi_{\mcT}(\mcP)\subset \mcT-\modu_{\infty}$ be the $n$-tilting subcategory of Thm. \ref{MiThm}, then $\mcC:=\kdual \widetilde{\mcT}\subset \modu_{\infty}-\mcT$ is an $n$-cotilting subcategory and we have 
\[ 
{}_{\perp} \widetilde{\mcT} = {}^{\perp} \mcC
\]
\end{cor}    
    
\subsection{Ideq tilting in relative homological algebra}
 We look at exact substructures with enough projectives on exact catgeories of the form $\modu_{\infty}-\mcP$ with $\mcP$ essentially small. Recall, if $\mcE=(\mcA, \mcS)$ is an exact category with underlying additive category $\mcA$ and class of short exact sequences $\mcS$, then an \textbf{exact substructure} is an exact category $\mcE'=(\mcA, \mcS')$ with $\mcS' \subset \mcS$. 
 
 We prove the following: 
 \begin{thm}
 Let $\mcP$ be an idempotent complete, 
 additive category. 
 Let $\mcE$ be an exact substructure of $\modu_{\infty}-\mcP$, with enough projectives $\mcQ:= \mcP(\mcE)$. Then, $\mcQ$ is $2$-ideq $0$-tilting subcategory of $\mcE$. 
 \end{thm}
  
 As a trivial corollary of Theorem \ref{ideqWep} and the previous Theorem, we obtain: 
 \begin{cor}\label{relIdeqWep}
  Let $\mcP$ be an idempotent complete, 
  additive category. 
 Let $\mcE$ be an exact substructure of $\modu_{\infty}-\mcP$, with enough projectives $\mcQ:= \mcP(\mcE)$. Then for every $n\geq 0$, every $n$-tilting subcategory of $\mcE$ is ideq $n$-tilting. 
 \end{cor}

 We prove two lemmata for the proof. 
 \begin{lem}
 In the previous situation. The functor $\mathbb{P}\colon \modu_{\infty}-\mcP \to \modu_{\infty}\mcQ, X\mapsto \Hom (-,X)\lvert_{\mcQ}$ has a left adjoint functor $\Phi'\colon \modu_{\infty}-\mcQ \to \modu_{\infty}-\mcP$ given by the restriction functor $\Phi' (X)= X\lvert_{\mcP}$. Furthermore, $\Phi'$ is exact.  
 \end{lem}
 
 \begin{proof}
 We consider $\mcQ \subset \modu_{\infty}- \mcP\subset \Modu-\mcP$. Then there is an adjoint pair of functors $\Phi \colon \Modu-\mcP \to \Modu- \mcQ\colon \Phi'$ with $\Phi (X)= \Hom_{\Modu- \mcP} (-,X)\lvert _{\mcQ} $ (cf. \cite{S}). By loc. cit. Cor.3.15, we have for $X \in \gen_{\infty}^{\modu_{\infty}-\mcP}(\mcQ ) = \modu_{\infty}-\mcP$ (by assumption) we have $\mathbb{P} (X)=\Phi (X) \in \modu_{\infty} \mcQ$ and $\Phi' (X) (P) = X(P)$. Therefore, the restriction functor is the left adjoint if it is well-defined.   We claim: If $ X\in \modu_{\infty}-\mcQ$, then $X \lvert_{\mcP}\in \modu_{\infty}-\mcP$.\\
 We apply the restriction functor to a projective resolution of $X$. This gives a right bounded complex $Q_*$ in $\modu_{\infty}-\mcP$ with terms in $\mcQ$ which is quasi-isomorphic to the restricted stalk complex of $X$. Now, by  \cite{Bue}, Thm 12.7, there exists a quasi-isomorphism $P_* \to Q_*$ with $P_*$ a right bounded complex of projectives in $\modu_{\infty}-\mcP$. Since compositions of quasi-isomorphisms are quasi-isomorphisms, the quasi-isomorphism  $P_*$ to the stalk complex $X\lvert_{\mcP}$ gives a projective resolution. 
 \end{proof}
 
 \begin{lem}
 Let $\mcE$ be an idempotent complete exact category with enough projectives given by $\mcQ$. If the functor 
 $\mathbb{P}\colon \mcE \to \modu_{\infty}-\mcQ, X \mapsto \Hom(-,X)\lvert_{\mcQ}$ has an exact left adjoint, then we have $\Reso_2(\Bild \mathbb{P})=\modu_{\infty}-\mcQ$. 
 \end{lem}
 
 \begin{proof}
 Let $X\in \modu_{\infty}-\mcQ$. We choose a projective presentation $\Hom_{\mcQ} (-, Q_1) \xrightarrow{\Hom(-,f)} \Hom_{\mcQ} (-, Q_0) \to X \to 0$ and denote $\Omega^2= \Ker (\Hom (-,f) )$.
 We claim that the unit $\Omega^2 \to \mathbb{P} \Phi' (\Omega^2)$ is an isomorphism (then $\Omega^2 \in \Bild (\mathbb{P})$ and since the projectives are in $\Bild \mathbb{P}$, the claim follows). \\
 First of all, we observe that any object $Z \in \Bild (\mathbb{P})$ fulfills that the unit is an isomorphism on $Z$ - this follows from the triangle identity of the adjunction. \\
 In particular, this holds for the projectives which lie in $\Bild (\mathbb{P})$.  
 Since $\mathbb{P} \circ \Phi'$ preserves kernels (since $\mathbb{P}$ preserves kernels and $\Phi'$ is exact) we can deduce from the commutative diagram 
 \[ 
 \xymatrix{
 0\ar[r]& \Omega^2 \ar[r]\ar[d]& \Hom_{\mcQ} (-, Q_1)\ar[r]\ar[d]& \Hom_{\mcQ} (-, Q_0 )\ar[d]\\
 0\ar[r]& \mathbb{P} \circ \Phi'(\Omega^2) \ar[r]& \mathbb{P} \circ \Phi'(\Hom_{\mcQ} (-, Q_1))\ar[r]&\mathbb{P} \circ \Phi' (\Hom_{\mcQ} (-, Q_0 ))
 }
 \]
 that the unit $\Omega^2 \to \mathbb{P} \circ \Phi' (\Omega^2)$ is an isomorphism. 
 \end{proof}

 \begin{exa}
 In the special case, $\mcA= \Lamod$ for an artin algebra $\La$, then 
 the definition of a \emph{relative tilting object} has been given in  \cite{ASoII}, the derived equivalence had been proven in \cite{Bu-Closed}. The proof in loc. cit. claims, we have an equivalence of additive categories and therefore, we have an equivalence on $K^{b,-}(-)$ of those - but this triangulated category depends not just on the additive category but on the ambient exact category and therefore a further explanation for the triangle equivalence (as given here, is helpful for the understanding). 
\end{exa}    

\section{Examples}

We plan to write an extended separate article on this topic. Therefore, we restrict to very view examples. This is my favourite construction of tilting subcategories: 

\begin{exa} \textbf{(special tilting)}
 
 \begin{lem} 
Let $\mcE$ be an exact category with enough projectives $\mcP$. Let $n\geq 1$ and $\mcM$ be full self-orthogonal subcategory, closed under summands with $\pdim_\mcE \mcM\leq 1$ and assume that \[\mcP\subset \cogen^{n-1} (\mcM). \] 
Let $\Omega^{-n}_{\mcM}\mcP$ denote the full subcategory of $\mcE$ consisting of objects $X$ such that there exists an exact sequence 
\[
0\to P \to M_0 \to \cdots \to M_{n-1} \to X \to 0 
\]
with $M_i \in \mcM, P \in \mcP$ and $\Hom_{\mcE}(-,\mcM)$ exact on it. 
We define $\mcT_n:= \mcM \vee \add(\Omega^{-n}_{\mcM}\mcP)$. 
Then $\mcT_n$ is an $n$-tilting subcategory and 
\[
\mcT_n^{\perp}= \gen_{n-1}(\mcM)
\]
\end{lem}
\begin{proof} The proof of the dual of Lem. 8.3 in \cite{MS} also carries through to show (t1), (t2) and (t3) (of Thm \ref{usualtilting}). 
\end{proof}
\begin{dfn}
Let $\mcE, \mcM$ be as in the previous lemma. If $\mcM\subset \mcP$, then we call $\mcT_n= \mcM \vee \add(\Omega^{-n}_{\mcM}\mcP)$ the $\mcM$-\textbf{special tilting subcategory}.
\end{dfn}

 \end{exa}

\begin{exa} \textbf{(Tilting modules for infinitely presented modules over rings)}
 Let $R$ be an associative unital ring. 
 We set $R-\Modu:=\Modu - (\rm{proj-R}),R-\modu_{\infty}:=\modu_\infty-(\rm{proj}-R)$ where $\rm{proj}-R$ denotes the category of finitely generated projective left $R$-modules. 
 This notation is justified by the observation that the category of all left $R$-modules is equivalent to $\Modu - (\rm{proj-R})$, just consider the following functor
 \[ 
 \begin{aligned}
 R-\Modu &\longrightarrow  \; \Modu - (\rm{proj-R})\\
 M &\mapsto \; \bigl(P \mapsto \Hom_R(P, M) \bigr)
 \end{aligned}
 \]
 It is an equivalence with quasi-inverse given by $F \mapsto F(R)$. 
Let now $\mcE:=R-\modu_{\infty}$. This is an exact category with enough projectives. 
 Let $T$ be an object in $\mcE$ and $\Gamma =\End_\mcE(T)^{op}$. 
We have $\add(T)$ has weak kernels if and only if $\Gamma$ is left coherent but we do not need to assume this here.\\
Then $T$ is $n$-tilting in $\mcE$ if and only if it satisfies (t1),(t2) and (t3) (cf. Theorem \ref{usualtilting}). By Thm \ref{ideqWep} we have an induced equivalence on bounded derived categories $D^b(R-\modu_{\infty})\to D^b(\Gamma-\modu_{\infty})$. This is also implied by Rickard's Morita theory for derived categories (\cite{R}, Thm 6.4 and Prop. 8.1).   \\

 
 
 \vspace{0.3cm}
\noindent 
On endomorphism rings of generators one can always find a special $1$-tilt:


Let $R$ be a ring and $M$ be a left $R$-module and $Q$ be a projective left $R$-module such that there is an epimorphism $Q^n\to M$ for some $n\geq 1$. Let $E= M\oplus Q$ and $\Gamma = \End_R(E)^{\op}$. Then, $P= \Hom_R(Q, E)$ is a projective right $\Gamma $-module. \\
Take the short exact sequence 
\[ 0\to K=\Ker (f) \to Q^{n+1} \to E \to 0\]
and apply the functor $\Hom_R(-, E)$ 
\[ 
0\to \Gamma \xrightarrow{F} P^{n+1} \to T_1:= \coKer (F) \to 0.
\]
We set $T=P \oplus T_1$. Then, it is straight forward to see:
 
\begin{cor}
$T\in \modu_{\infty}-\Gamma$ is a the special $1$-tilting module for $\mcM= \add(P)$.  In particular, ${\gen}(T) = {\gen}(P)$.    
\end{cor}
 \end{exa}


\bibliographystyle{alpha}
\bibliography{Sauter-TiltExact}
\end{document}